\documentclass[12pt,oneside,reqno]{amsart}

\usepackage{amssymb}
\usepackage{amsfonts}
\usepackage{bbm}
\usepackage{cases}
\usepackage{amsmath}
\usepackage{graphicx}
\usepackage{mathrsfs}
\usepackage{stmaryrd}
\usepackage{txfonts}
\usepackage{enumerate,amsmath,amssymb,amsthm,booktabs}

\numberwithin{equation}{section}

\newcommand{\be}{\begin{eqnarray}}
\newcommand{\ee}{\end{eqnarray}}
\newcommand{\ce}{\begin{eqnarray*}}
\newcommand{\de}{\end{eqnarray*}}
\newtheorem{theorem}{Theorem}[section]
\newtheorem{lemma}[theorem]{Lemma}
\newtheorem{remark}[theorem]{Remark}
\newtheorem{definition}[theorem]{Definition}
\newtheorem{proposition}[theorem]{Proposition}
\newtheorem{Examples}[theorem]{Example}
\newtheorem{corollary}[theorem]{Corollary}

\def\Re{{\mathrm{Re}}}

\def\eps{\varepsilon}

\def\e{\mathrm{e}}

\def\p{\partial}

\def\[{{\Big[}}
\def\]{{\Big]}}
\def\<{{\langle}}
\def\>{{\rangle}}
\def\({{\Big(}}
\def\){{\Big)}}

\def\bx{{\mathbf{x}}}
\def\tr{\mathrm {tr}}

\def\dif{{\mathord{{\rm d}}}}

\def\={&\!\!=\!\!&}

\def\cM{{\mathcal M}}

\def\cS{{\mathcal S}}
\def\cT{{\mathcal T}}

\def\mE{{\mathbb E}}

\def\mI{{\mathbb I}}

\def\mM{{\mathbb M}}
\def\mN{{\mathbb N}}

\def\mQ{{\mathbb Q}}
\def\mR{{\mathbb R}}

\def\1{{\mathbf{1}}}

\def\sL{{\mathscr L}}

\def\sP{{\mathscr P}}
\def\sQ{{\mathscr Q}}

\def\wt{\widetilde}

\def\geq{\geqslant}
\def\leq{\leqslant}

\def\Re{{\mathrm{Re}}}

\def\eps{\varepsilon}

\def\e{\mathrm{e}}

\def\p{\partial}

\def\[{{\Big[}}
\def\]{{\Big]}}
\def\<{{\langle}}
\def\>{{\rangle}}
\def\({{\Big(}}
\def\){{\Big)}}

\def\bx{{\mathbf{x}}}
\def\tr{\mathrm {tr}}

\def\dif{{\mathord{{\rm d}}}}

\def\={&\!\!=\!\!&}
\def\bt{\begin{theorem}}
\def\et{\end{theorem}}
\def\bl{\begin{lemma}}
\def\el{\end{lemma}}
\def\br{\begin{remark}}
\def\er{\end{remark}}
\def\bx{\begin{Examples}}
\def\ex{\end{Examples}}
\def\bd{\begin{definition}}
\def\ed{\end{definition}}
\def\bp{\begin{proposition}}
\def\ep{\end{proposition}}
\def\bc{\begin{corollary}}
\def\ec{\end{corollary}}

\def\geq{\geqslant}
\def\leq{\leqslant}

\def\bX{{\mathbf X}}

 \def\R{\mathbb R}
 \def\R{\mathbb R}

\def\<{\langle} \def\>{\rangle}

 \def\beq{\begin{equation}}  
 
\def\e{\text{\rm{e}}}

\allowdisplaybreaks

\begin{document}

\title{Propagation of regularity in $L^p$-spaces for Kolmogorov  type hypoelliptic operators}

\date{}

\author{{Zhen-Qing Chen} \ \
and \ \ {Xicheng Zhang}}

\address{Zhen-Qing Chen:
Department of Mathematics, University of Washington, Seattle, WA 98195, USA\\
Email: zqchen@uw.edu
 }

\address{Xicheng Zhang:
School of Mathematics and Statistics, Wuhan University,
Wuhan, Hubei 430072, P.R.China\\
Email: XichengZhang@gmail.com
 }

\begin{abstract}
Consider the following   
Kolmogorov  type hypoelliptic operator
$$
\sL_t := \mbox{$\sum_{j=2}^n$}x_j\cdot\nabla_{x_{j-1}}+\tr(a_t \cdot\nabla^2_{x_n})
 $$
on $\mR^{nd}$, where $n\geq 2$, $d\geq 1$, 
  $x=(x_1,\cdots,x_n)\in (\mR^d)^n =\mR^{nd}$ and $a_t$ is a time-dependent constant
  symmetric  $d\times d$-matrix that is uniformly elliptic and bounded.
Let $\{\cT_{s, t};  t\geq s\}$ be the time-dependent semigroup associated with $\sL_t$;  that is,
$\partial_s  \cT_{s, t}  f = - \sL_s \cT_{s, t}f$.   
For any $p\in(1,\infty)$, we show that there is a constant $C=C(p,n,d)>0$ 
such that for any $f(t, x)\in L^p(\mR \times \mR^{nd})=L^p(\mR^{1+nd})$
and every $\lambda \geq 0$, 
$$
 \left\|\Delta_{x_j}^{ {1}/ {(1+2(n-j))}}  \int^{\infty}_0 \mathrm{e}^{-\lambda t} \cT_{s, t+s}f(t+s, x)\dif t\right\|_p\leq C\|f\|_p,\quad  j=1,\cdots, n,
$$
where $\|\cdot\|_p$ is the usual $L^p$-norm in $L^p(\mR\times \mR^{nd}; \dif s\times\dif x)$. To show this type of estimates, we first study the propagation of
regularity in $L^2$-space from variable $x_n$ to $x_j$, $1\leq j\leq n-1$,  
 for the solution of  the transport equation  
$ \p_t u+ \sum_{j=2}^nx_j\cdot\nabla_{x_{j-1}} u=f.$

\bigskip
\noindent 
\textbf{Keywords}: 
Kolmogorov's hypoelliptic operators, Fefferman-Stein's theorem, Propagation of $L^p$-regularity\\

\noindent
 {\bf AMS 2010 Mathematics Subject Classification:}  Primary: 42B20, 60H30, Secondary 35H10, 35Q20
 
 \end{abstract}

\thanks{
This work is partially supported by Simons Foundation grant 520542 and by NNSFC grant of China (No. 11731009). }

\maketitle

\section{\bf Introduction}

Let $n\geq 2$ and $d\in\mN$. 
Denote by  $\mM^d_{sym}$   the set
of all symmetric $d\times d$-matrices. 
In this paper we consider the following Kolmogorov type hypoelliptic operator
on $\mR^{nd}$:
\begin{align}\label{SL}
\sL_t:=\sum_{i,j=1}^d a^{ij}_t\p_{x_{ni}}\p_{x_{nj}}+\sum_{j=2}^n x_j\cdot\nabla_{x_{j-1}},
\end{align}
where $x=(x_1,x_2,\cdots,x_n)\in\mR^{nd}$ with $x_j=(x_{j1},\cdots,x_{jd})\in\mR^d$ for each $j=1,\cdots,n$,
$\nabla_{x_j}=(\p_{x_{j1}},\cdots,\p_{x_{jd}})$,
and $a_t=(a^{ij}_t):\mR\to \mM^d_{sym}$ is a measurable map
(independent of $x$)  having
\begin{align}\label{Ell}
\kappa^{-1}\mI_{d\times d}\leq a_t\leq\kappa\mI_{d\times d}
\end{align}
for some $\kappa\geq 1$.
Here $\mI_{d\times d}$ stands for the $d\times d$ identity matrix.
Let $\nabla:=(\nabla_{x_1},\cdots,\nabla_{x_n})$, $\nabla^2_{x_n}:=(\p_{x_{ni}}\p_{x_{nj}})_{i,j=1,\cdots,d}$ and
\begin{align}\label{AA}
A=A_n=\left(
\begin{array}{ccccc}
0_{d\times d}&\mI_{d\times d}&0_{d\times d}&\cdots&\cdots  \smallskip \\
0_{d\times d}&0_{d\times d}&\mI_{d\times d}&0_{d\times d}&\cdots \smallskip \\
\vdots&\ddots&\ddots&\ddots&\vdots  \smallskip \\
\vdots&\cdots&\ddots&0_{d\times d}&\mI_{d\times d}  \smallskip \\
0_{d\times d}&\cdots&\cdots&0_{d\times d}&0_{d\times d}
\end{array}
\right)_{nd\times nd}.
\end{align}
We can rewrite $\sL_t$ as the following compact form:
$$
\sL_t=\tr(a_t\cdot\nabla^2_{x_n})+Ax\cdot\nabla,
$$
where ``tr'' denotes the trace of a matrix. 
The hypoelliptic operator   
$\sL_t$ is a differential operator 
of second order in $x_n$  variable but is of order 1 in variables $x_1, \cdots, x_{n-1}$.
These variables are connected through the drift terms $x_j\cdot\nabla_{x_{j-1}}$ for $2\leq j\leq n$. 
In this paper, we study the regularity of the resolvent functions of $\sL$.
Roughly speaking, we show that Laplacian in $x_n$-variable of the resolvent of $\sL$ is a bounded operator in $L^p$-space, and that
this regularity in $x_n$ propagates  to other variables  in such a way that   the fractional Laplacian in $x_j$-variable of power  $1/(1+2(n-j))$  
of the resolvent of $\sL$ 
is a bounded operator in $L^p$-space  for $1\leq j \leq n-1$; see Theorem \ref{Main} below for a precise statement.

\medskip

Consider the following linear stochastic differential equations (SDEs)
in $\R^{nd}$:
\begin{align}\label{SDE}
\dif X^{s,x}_t =AX^{s,x}_t \dif t+\sigma^a_t\dif W_t \quad \hbox{for } t>s 
 \ \hbox{ with }  X^{s,x}_s=x \in (\R^d)^n, 
\end{align}
where $(W_t)_{t\in\mR}$ is a standard $nd$-dimensional Brownian motion and
\begin{align}\label{ES11}
\sigma_t^a:=\left(
\begin{array}{cc}
0_{(n-1)d\times(n-1)d},&0_{(n-1)d\times d} \smallskip \\
0_{d\times(n-1)d},&(\sqrt{2a_t})_{d\times d}
\end{array}
\right)_{nd\times nd}.
\end{align}
It is easy to see that the solution $X^{s,x}_t$ of \eqref{SDE} is explicitly given by 
$$
X^{s,x}_t = \e^{(t-s)A}x+\int^t_s\e^{(t-r)A}\sigma^a_r\dif W_r,
 $$
where $\e^{tA}$ is the exponential matrix of $A$ with the expression
\begin{align}\label{ES1}
\e^{tA}=\left(
\begin{array}{ccccc}
\mI_{d\times d}&t\mI_{d\times d}&\tfrac{t^2}{2} \mI_{d\times d} &\cdots&\tfrac{t^{n-1}}{(n-1)!} \mI_{d\times d} \smallskip \\
0_{d\times d}&\mI_{d\times d}&t\mI_{d\times d}&\cdots&\tfrac{t^{n-2}}{(n-2)!} \mI_{d\times d} \smallskip \\
\vdots&\ddots&\ddots&\ddots&\vdots \smallskip \\
\vdots&\cdots&0_{d\times d}&\mI_{d\times d}&t\mI_{d\times d}  \smallskip \\
0_{d\times d}&\cdots&\cdots&0_{d\times d}&\mI_{d\times d}
\end{array}
\right)_{nd\times nd}.
\end{align}
Notice that if $a_t=a$ does not depend on $t$ (i.e., time homogeneous), then 
\begin{equation}\label{e:1.7}
X^{s,x}_t\stackrel{(d)}{=}Z^x_{t-s} , 
\quad  \hbox{where   }  \  Z^x_t:=
 \e^{tA}x+\int^{t}_0\e^{rA}\sigma^a\dif W_r.
\end{equation} 
In this case, $Z^x_t$ is an $(nd)$-dimensional Gaussian random variable with density
$$
y\mapsto p_t(x,y)=\frac{\e^{-(\Theta_{t^{-1/2}}(y-\e^{tA}x))^*\Sigma^{-1} \Theta_{t^{-1/2}}(y-\e^{tA}x)}}{((2\pi)^{nd}t^{n^2d}\det(\Sigma))^{1/2}},
$$
where $\Theta_r:\mR^{nd}\to\mR^{nd}$ is the dilation operator defined by
\begin{align}\label{Dia}
\Theta_r(x)=(r^{2n-1}x_1,r^{2n-3}x_2,\cdots, rx_n),
\end{align}
and $\Sigma:=\int^1_0\e^{r A}\sigma^a(\sigma^a)^*\e^{rA^*}\dif r$ is the covariance matrix of $Z$  (see \cite{Ko} or \eqref{HJ1} below).

For $f\in C^2_b(\mR^{nd})$, define
\begin{align}\label{TT}
\cT_{s,t}f(x):=\mE f(X^{s,x}_t).
\end{align}
It is easy to check that
$\sL_s$ 
 is the infinitesimal generator of $\cT_{s,t}$, that is, 
\begin{align}\label{Gen}
\p_s\cT_{s,t}f+\sL_s\cT_{s,t}f=0.
\end{align}

The goal of this paper   is to
establish  the following $L^p$-maximal regularity estimate.
 
\bt\label{Main}
 Let $p\in(1,\infty)$. 
Under the uniform ellipticity condition \eqref{Ell},
 there is a constant $C=C(n,\kappa,p,d)>0$ such that for any 
 $f\in L^p(\mR\times \mR^{nd})$
 and $\lambda\geq 0$,
\begin{align}\label{BN5}
\left\|\Delta_{x_j}^{ {1}/{(1+2(n-j))}}
 \int^{\infty}_0 \e^{ - \lambda t}\cT_{s,t+s}f(t+s,x )
\dif t\right\|_p\leq C\|f\|_p,
\quad j=1, \dots, n, 
\end{align}
where  $\Delta_{x_j}^{ {1}/{(1+2(n-j))}}:=-(-\Delta_{x_j})^{{1}/{(1+2(n-j))}}$ is the fractional Laplacian acting 
on the $j$-th variable $x_j\in \mR^d$. 
\et
Note that 
 $$
u(s, x):=  \int^{\infty}_0 \e^{ - \lambda t}\cT_{s,t+s}f(t+s,x ) \dif t 
= \int_s^\infty  \e^{ - \lambda (r-s)} \cT_{s,r}f(r,x ) \dif r 
$$
satisfies 
\begin{equation}\label{e:1.12} 
 \partial_s u (s, x) + (\sL_s-\lambda ) u (s, x) + f(s, x)=0.
 \end{equation}
One of the motivation of studying the estimate \eqref{BN5} comes from the study of the following $(n+1)$-order stochastic differential equation:
\begin{align}\label{SDE0}
\dif X^{(n)}_t=b_t(X_t,X^{(1)}_t,\cdots, X^{(n)}_t)\dif t+\sigma_t(X_t,X^{(1)}_t,\cdots, X^{(n)}_t)\dif \wt W_t , 
\end{align}
where $X^{(n)}_t$ denotes the $n$th-order derivative of $X_t$ with respect to the time variable, 
$b:\mR_+\times\mR^{(n+1)d}\to\mR^d$ and $\sigma: \mR_+\times\mR^{(n+1)d}\to\mR^d\otimes\mR^d$ are  measurable functions,
and $\wt W_t$ is a $d$-dimensional Brownian motion.  
Notice that if we let
$$
\bX_t:=(X_t, X^{(1)}_t,\cdots, X^{(n)}_t),
$$
then $\bX_t$ solves the following one order stochastic differential equation
$$
\dif \bX_t=(X^{(1)}_t,\cdots, X^{(n)}_t, b_t(\bX_t))\dif t+(0,\cdots, 0,\sigma_t(\bX_t)\dif \wt W_t ), 
\quad \mathbf X_0=\mathbf x,
$$
where $\mathbf x=(\mathbf x_{i})_{i=0,\cdots,n}=((x_{ij})_{j=1,\cdots, d})_{i=0,\cdots,n}$.
In particular, the infinitesimal generator of Markov process $\bX_t(\mathbf x)$ is given by
$$
\sL_t f(\mathbf x)=\sum_{i,j,k=1}^d (\sigma^{ik}_t\sigma^{jk}_t)(\mathbf x)\p_{x_{ni}}\p_{x_{nj}}f(\mathbf x)
+\sum_{j=1}^n\mathbf x_j\cdot\nabla_{\mathbf x_{j-1}} f(\mathbf x)+b_t(\mathbf x)\cdot\nabla_{\mathbf x_n} f(\mathbf x),
$$ 
which is of the form  similar to \eqref{SL}. 
Thus, the estimate \eqref{BN5} could be used to study the well-posedness of SDE \eqref{SDE0}
with rough coefficients $b$ and $\sigma$. Indeed, when $n=1$ and $\sigma$ is bounded and uniformly nondegenerate, 
the second named author \cite{Zh2} 
studied the strong well-posedness of SDE \eqref{SDE0}
with both $(\mI  -\Delta_{\bf x_1})^{1/3}b$ and $\nabla\sigma$ in $L^p_{loc}(\mR_+\times\mR^{2d})$
for some $p>4d+2$. See also 
\cite{Fe-Fl-Pr-Vo} for similar results when $\sigma_t=\mI_{d\times d}$.

We  now recall some related 
results in literature 
about the estimate \eqref{BN5}. In \cite{Br-Cu-La-Pr}, 
the authors adopted Coifman-Weiss' theorem to show the estimate \eqref{BN5} for $j=n$.
When $n=2$, in \cite{Ch-Zh}, we established a version of  Fefferman-Stein's theorem and then used it 
to show  the estimate \eqref{BN5} for $j=1,2$ even for nonlocal operators. 
It should be noticed that the methods used in \cite{Br-Cu-La-Pr} and \cite{Ch-Zh} are quite different. 
In \cite{Br-Cu-La-Pr}, the key point is to show some weak 1-1 type estimate. While in \cite{Ch-Zh},
the main point is to show that the operator in \eqref{BN5} is bounded from $L^\infty$ to some $BMO$ spaces.
In particular, to show the propagation of the regularity from the nondegenerate component to the degenerate component, in \cite{Ch-Zh}, we have used Bouchet's result \cite{Bo}. More precisely, Bouchet studied the following transport equation:
$$
\p_t u+x_2\cdot\nabla_{x_1} u=f,
$$
and showed that for any $\alpha\geq 0$,
$$
\|\Delta_{x_1}^{\frac{\alpha}{2(1+\alpha)}}u\|_2\leq C\|\Delta_{x_{2}}^{\frac{\alpha}{2}}u\|^{\frac{1}{1+\alpha}}_2\|f\|_2^{\frac{\alpha}{1+\alpha}},
$$
where $C=C(\alpha,d)>0$.
A simplified proof of this type estimate was provided in \cite{Al}. Thus, the first 
goal of this paper is to extend  the above estimate to the following more general transport equation:
$$
\p_t u+\sum_{j=2}^nx_j\cdot\nabla_{x_{j-1}} u=f.
$$
That is, we want to show that for any $j=1,\cdots, n-1$ and $\alpha\geq 0$, there is a constant $C=C(\alpha,d,j,n)>0$ such that 
$$
\|\Delta_{x_j}^{\frac{\alpha}{2(1+\alpha)}}u\|_2\leq C\|\Delta_{x_{j+1}}^{\frac{\alpha}{2}}u\|^{\frac{1}{1+\alpha}}_2\|f\|_2^{\frac{\alpha}{1+\alpha}}.
$$
Such an extension from $n=2$ to $n\geq 3$ is non-trivial, 
see Section 3.

Although the above result is proven for Laplacian operator, we can extend it 
 to  more general nonlocal operator as in \cite{Ch-Zh} without   any difficulty. 
Indeed, consider the following nonlocal operator:
$$
\widetilde\sL^\nu_\sigma f(x):=\int_{\mR^d} [f(x+\sigma y)+f(x-y)-2f(x)]\nu(\dif y),
$$
where $\sigma\in\mM^d$ is a $d\times d$ matrix and $\nu$ is a symmetric L\'evy
measure on $\mR^d$ (that is, $\nu$ is a measure on $\mR^d \setminus \{0\}$ with 
$\nu (A)=\nu (-A)$ and $\int_{\mR^d \setminus \{0\}} (1\wedge |z|^2) \nu (\dif z)<\infty$). 
Let $n\geq 2$ and 
$$
\sL_tf(x):=\widetilde\sL^{\nu_t}_{\sigma_t, x_n}f(x)+\sum_{j=2}^nx_j\cdot\nabla_{x_{j-1}} f(x),
$$
where $\widetilde\sL^{\nu_t}_{\sigma_t, x_n}$ means that the operator acts on the variable $x_n$.
Suppose that 
$$
\|\sigma\|_\infty+\|\sigma^{-1}\|_\infty<\infty
$$
and for some $\alpha\in(0,2)$,
$$
\nu^{(\alpha)}_1\leq\nu_s\leq\nu^{(\alpha)}_2,
$$
where $\nu^{(\alpha)}_1$ and $\nu^{(\alpha)}_2$ are two symmetric and nondegenerate $\alpha$-stable L\'evy measures (see \cite{Ch-Zh}).
Under the above assumptions, 
we can show as in \cite{Ch-Zh}  that for any $j=1,\cdots, n$,
\begin{align}\label{BN05}
\left\|\Delta_{x_j}^{\frac{\alpha}{2(1+\alpha(n-j))}}\int^{\infty}_0\e^{-\lambda t}\cT^{\nu,\sigma}_{s,t+s}f(t+s,x)\dif t\right\|_p\leq C\|f\|_p,
\end{align}
where $\cT^{\nu,\sigma}_{s,t}$ is defined as in \eqref{TT} by using 
the time-inhomogeneous Markov process 
$\{\{Z^{s, x}_t; t\geq 0\}; (s, x)\in \mR\times \mR^{nd}\}$
determined by the family of L\'evy measures $\{\nu_s, s\in \mR\}$ in place of  Brownian motion.
We note 
 that at the almost same time, Huang, Menozzi and Priola \cite{Hu-Me-Pr}
 have obtained  
\eqref{BN05} for time-independent $\sigma$ and $\nu$ by using Coifman-Weiss' theorem. 
As mentioned above, 
our proof is based on 
a tailored version of Fefferman-Stein's theorem. 

This paper is organized as follows: In Section 2, we prepare some estimates about 
the probability density function of $X^{s,x}_t$, and establish a Fefferman-Stein type  theorem.
In Section 3, we show the propagation of the regularity for transport equation. In Section 4, we prove our main result.

Throughout this paper, we use the following convention: The letters $C$ and $c$ with or without subscripts will denote a positive constant, 
whose value may change in different places.
Moreover, we use $A\lesssim B$ to denote $A\leq CB$ for some constant $C>0$.

\section{\bf Preliminaries}

\subsection{Estimate of density of  $X^{s,x}_{t}$}

In this subsection we establish some estimates on the density of  $X^{s,x}_{t}$,  
which will be used later. 

\bl\label{Th23}
Under \eqref{Ell}, 
$X^{s, {\bf 0}}_t$ of \eqref{e:1.7}  has a smooth density function $p^{({\bf 0})}_{s,t}(y)$.  For
each $\beta=(\beta_1,\cdots,\beta_n)\in\mN_0^n$,
where $\mN_0=\{0\}\cup\mN$,
there are constants $C, c>0$ only depending on $n,\beta, d$ and $\kappa$ such that for all $s<t$ and  
$y\in\mR^{nd}$,
\begin{align}\label{HJ4}
|\nabla^{\beta_1}_{y_1}\cdots\nabla^{\beta_n}_{y_n}p^{({\bf 0})}_{s,t} (y)|\leq C(t-s)^{-(n^2d+\sum_{i=1}^n(2(n-i)+1)\beta_i)/2}\e^{-c|\Theta_{(t-s)^{-1/2}}y|^2},
\end{align}
where $\Theta_r$ is the dilation operator defined by \eqref{Dia}.
\el

\begin{proof}
Since $(cW_{\cdot/c^2})\stackrel{(d)}{=}W$ for $c\not=0$, by a  change of variables, we have
$$ 
X^{s, {\bf 0}}_t 
=\int^t_s\e^{(t-r)A}\sigma^a_r\dif W_r\stackrel{(d)}{=}(t-s)^{1/2}\int^1_0\e^{(t-s)(1-r)A}\sigma^a_{s+(t-s)r}\dif W_{r}.
$$
Hence, by definitions \eqref{ES11}, \eqref{ES1} and \eqref{Dia},
$$
\Theta_{(t-s)^{-1/2}}
X^{s, {\bf 0}}_t 
\stackrel{(d)}{=}\int^1_0\e^{(1-r)A}\sigma^a_{s+(t-s)r}\dif W_r
\stackrel{(d)}{=}\int^1_0\e^{rA}\tilde\sigma^a_r\dif W_r=:Z,
$$
where $\tilde\sigma^a_r:=\sigma^a_{s+(t-s)(1-r)}$. Since $Z$ is a $nd$-dimensional Gaussian random variable with mean value zero and covariance matrix
$$
\Sigma=\int^1_0\e^{rA}\tilde\sigma^a_r(\tilde\sigma^a_r)^*\e^{rA^*}\dif r,
$$
we have
\begin{align}\label{HJ1}
p^{({\bf 0})}_{s,t}(y) 
=\frac{\e^{-(\Theta_{(t-s)^{-1/2}}y)^*\Sigma^{-1} (\Theta_{(t-s)^{-1/2}}y)}}{((2\pi)^{nd}(t-s)^{n^2d}\det(\Sigma))^{1/2}}.
\end{align}
On the other hand, by \eqref{ES11}, \eqref{ES1} and \eqref{AA}, we have for all $y\in\mR^{nd}$,
\begin{align*}
\begin{split}
y^*\Sigma y&=\int^1_0|y^*\e^{rA}\tilde\sigma^a_r|^2\dif r\geq2\kappa^{-1}\int^1_0\left|\sum_{j=1}^n \frac{r^{n-j}}{(n-j)!}y_j\right|^2\dif r\\
&\geq2\kappa^{-1}|y|^2\inf_{|\omega|=1}\int^1_0\left|\sum_{j=1}^n \frac{r^{n-j}}{(n-j)!}\omega_j\right|^2\dif r.
\end{split}
\end{align*}
Since the unit sphere in $\mR^{nd}$ is compact, and for each $\omega \in \mR^{nd} \setminus \{{\bf 0}\}$, 
$$
\delta(\omega):=\int^1_0\left|\sum_{j=1}^n \frac{r^{n-j}}{(n-j)!}\omega_j\right|^2\dif r>0,
$$
 we have $c_0:=\inf_{|\omega|=1}\delta(\omega) >0$, and so 
\begin{align}\label{HJ2}
y^*\Sigma y\geq 2\kappa^{-1}|y|^2\inf_{|\omega|=1}\delta(\omega) \geq  2  \kappa^{-1} c_0 \, |y|^2.
\end{align}
The desired estimate now follows by the chain rule, \eqref{HJ1} and \eqref{HJ2}.
\end{proof}
For $\alpha\in(0,2]$, the fractional Laplacian $\Delta^{\alpha/2}$ in $\mR^d$ is defined by 
Fourier's transform as
$$
\widehat{\Delta^{\alpha/2} f}(\xi)=- |\xi|^\alpha\hat f(\xi),\ \ f\in\cS(\mR^d), 
$$
where $\cS(\mR^d)$ is the space of Schwartz rapidly decreasing functions. For $\alpha\in(0,2)$,
up to a multiplying constant, an alternative definition of $\Delta^{\alpha/2}$ is given by the following integral form (cf. \cite{St}):
\begin{align}\label{Loc}
\Delta^{\alpha/2} f(x)=\int_{\mR^d}\frac{\delta^{(2)}_zf(x)}{|z|^{d+\alpha}}\dif z,
\end{align}
where $\delta^{(2)}_zf(x):=f(x+z)+f(x-z)-2f(x).$ Observe that for $\alpha\in(0,1)$,
\begin{align}\label{Loc1}
\Delta^{\alpha/2} f(x)=2\int_{\mR^d}\frac{\delta^{(1)}_zf(x)}{|z|^{d+\alpha}}\dif z\mbox{ with }\delta^{(1)}_zf(x):=f(x+z)-f(x).
\end{align}
\bc\label{Cor1}
For any $j=1,\cdots,n$, $\alpha\in(0,2]$ and $\beta=(\beta_1,\cdots,\beta_n)\in\mN_0^n$, 
there is a constant $C>0$ such that for all $f\in C_b^\infty(\mR^{nd})$ and $s<t$,
\begin{align}
&\|\Delta^{\alpha/2}_{x_j}\nabla^{\beta_1}_{x_1}\cdots\nabla^{\beta_n}_{x_n}\cT_{s,t}f\|_\infty
\leq C|t-s|^{-(\sum_{i=1}^n(2(n-i)+1)\beta_i+(2(n-j)+1)\alpha)/2}\|f\|_\infty,\label{EK2}\\
&\|\nabla^{\beta_1}_{x_1}\cdots\nabla^{\beta_n}_{x_n}\cT_{s,t}\Delta^{\alpha/2}_{x_j}f\|_\infty
\leq C|t-s|^{-(\sum_{i=1}^n(2(n-i)+1)\beta_i+(2(n-j)+1)\alpha)/2}\|f\|_\infty,\label{EK22}
\end{align}
where $\Delta^{\alpha/2}_{x_j}$ means that the fractional Laplacian acts on the variable $x_j$.
\ec
\begin{proof}
Below we only show \eqref{EK2} and \eqref{EK22} for $\alpha\in(0,2)$. 
Let $p_{s,t}(x,y)$ be the probability density function 
of $X^{s,x}_t=\e^{(t-s)A}x+X^{s,{\bf 0}}_t$.
We have
\begin{align}\label{H1}
p_{s,t}(x,y)=p^{({\bf 0})}_{s,t}(y-\e^{(t-s)A}x).
\end{align}
For notational convenience, we write 
$$
h_y(x):=\nabla^{\beta_1}_{x_1}\cdots\nabla^{\beta_n}_{x_n}p_{s,t}(x,y),\ \gamma:=\mbox{$\sum_{i=1}^n$}(2(n-i)+1)\beta_i/2
$$
and 
$$
\delta^{(2)}_{z_j}h_y(x)=h_y(x+\tilde z_j)+h_y(x-\tilde z_j)-2h_y(x),\ \ \tilde z_j=(0,\cdots, z_j,\cdots, 0).
$$
By \eqref{HJ4}, \eqref{H1} and the chain rule, we have
\begin{align}\label{RR1}
\begin{split}
|\delta^{(2)}_{z_j}h_y(x)|&\leq C(t-s)^{-n^2d/2-\gamma}\Big(\e^{-c|\Theta_{(t-s)^{-1/2}}(y-\e^{(t-s)A}(x+\tilde z_j))|^2}\\
&+\e^{-c|\Theta_{(t-s)^{-1/2}}(y-\e^{(t-s)A}(x-\tilde z_j))|^2}+\e^{-c|\Theta_{(t-s)^{-1/2}}(y-\e^{(t-s)A}x)|^2}\Big),
\end{split}
\end{align}
and also by the mean value formula, 
\begin{align}\label{RR2}
|\delta^{(2)}_{z_j}h_y(x)|\leq C(t-s)^{-n^2d/2-\gamma-(2(n-j)+1)}\e^{-c|\Theta_{(t-s)^{-1/2}}(y-\e^{(t-s)A}\tilde x)|^2}|z_j|^2
\end{align}
for some $\tilde x\in\mR^{nd}$ depending on $z_j$. By formula  \eqref{Loc}, we have
$$
\Delta^{\alpha/2}_{x_j}h_y(x)
=\left(\int_{|z_j|>(t-s)^{(2(n-j)+1)/2}}+\int_{|z_j|\leq(t-s)^{(2(n-j)+1)/2}}\right)\frac{\delta^{(2)}_{z_j}h_y(x)}{|z_j|^{d+\alpha}}\dif z_j
=:I_1(x,y)+I_2(x,y).
$$
For $I_1(x,y)$, we have by \eqref{RR1}  
\begin{align*}
\int_{\mR^{nd}}|I_1(x,y)|\dif y&\lesssim(t-s)^{-\gamma}\int_{|z_j|>(t-s)^{(2(n-j)+1)/2}}\frac{\dif z_j}{|z_j|^{d+\alpha}}
\lesssim(t-s)^{-\gamma-(2(n-j)+1)\alpha/2},
\end{align*}
where we have used that 
$$
(t-s)^{-n^2d/2}\int_{\mR^{nd}}\e^{-c|\Theta_{(t-s)^{-1/2}}(y)|^2}\dif y=\int_{\mR^{nd}}\e^{-c|y|^2}\dif y.
$$
For $I_2(x,y)$, we have by \eqref{RR2}  
\begin{align*}
\int_{\mR^{nd}}|I_2(x,y)|\dif y&\lesssim(t-s)^{-\gamma-(2(n-j)+1)}\int_{|z_j|\leq(t-s)^{(2(n-j)+1)/2}}\frac{|z_j|^2}{|z_j|^{d+\alpha}}\dif z_j\\
&\lesssim(t-s)^{-\gamma-(2(n-j)+1)\alpha/2}.
\end{align*}
Combining the above calculations, we obtain
$$
|\Delta^{\alpha/2}_{x_j}\nabla^{\beta_1}_{x_1}\cdots\nabla^{\beta_n}_{x_n}\cT_{s,t}f(x)|
=\left|\int_{\mR^d}\Delta^{\alpha/2}_{x_j}h_y(x) f(y)\dif y\right|\leq C(t-s)^{-\gamma-(2(n-j)+1)\alpha/2}\|f\|_\infty.
$$
Thus we proved \eqref{EK2}. Similarly, we can show \eqref{EK2}.
\end{proof}

\subsection{Fefferman-Stein's theorem}
First of all, we introduce a family of ``balls'' in $\mR^{1+nd}$
that matches the geometry induced by the hypoelliptic operator \eqref{SL}. 
 For any $r>0$ and point $(t_0, x_0)\in\mR^{1+nd}$, we define
\begin{align*}
Q_r(t_0,x_0)&:=\Big\{(t,x): \ell(t-t_0,x-\e^{(t-t_0)A}x_0)\leq r\Big\},
\end{align*}
where 
\begin{align*}
\ell(t,x)&:=\max\left\{  | t|^{ {1}/{2}}, |x_1|^{{1}/{(2n-1)}},|x_2|^{{1}/{(2n-3)}},\cdots, |x_{n-1}|^{{1}/{3}},|x_n|\right\}.
\end{align*}
We use  $\mQ$ to denote the  set of all such balls. 

\bl\label{Le31}
\begin{enumerate}[\rm (i)]
\item $\ell(r^2t, \Theta_rx)=r\ell(t,x)$ for any $r>0$, where $\Theta_r$ is the dilation operator defined by
\eqref{Dia}.
\item $|Q_r(t_0,x_0)|=\omega^n_{d}r^{n^2d+2}$, where $|\cdot|$ denotes the Euclidean volume and $\omega_d$ is the volume of the unit ball in 
$(\mR^d, | \cdot |)$.  
\item For all $(t,x),(s,y),(r,z)\in\mR^{1+nd}$, we have
\begin{align}\label{Tri}
\begin{split}
&\ell(s-t,y-\e^{(s-t)A}x)\leq 3\ell(t-s,x-\e^{(t-s)A}y)\leq\\
&\quad\leq 12\Big(\ell(t-r,x-\e^{(t-r)A}z)+\ell(r-s,z-\e^{(r-s)A}y)\Big).
\end{split}
\end{align}
\item Suppose that $Q_r(t_0,x_0)\cap Q_r(t_0',x_0')\not=\emptyset$, then
\begin{align}
Q_r(t_0,x_0)\subset Q_{20\cdot r}(t_0',x_0').\label{EL1}
\end{align}
\item For $(t,x), (s,y)\in\mR^{1+nd}$, define
$$
\rho((t,x),(s,y)):=\ell(t-s,x-\e^{(t-s)A}y)+\ell(s-t,y-\e^{(s-t)A}x),
$$
and for $(t_0,x_0)\in\mR^{1+nd}$ and $r>0$,
$$
\tilde Q_r(t_0,x_0):=\{(t,x): \rho((t,x),(t_0,x_0))\leq r\}.
$$
Then $Q_{r}(t_0,x_0)\subset \tilde Q_{r}(t_0,x_0)\subset Q_{4r}(t_0,x_0)$.
\end{enumerate}
\el

\begin{proof}
(i) and (ii) follow directly from the definition of $\ell (t, x)$. 
\medskip\\
(iii) We only prove the second inequality in \eqref{Tri}. The first one is similar. Observing that for all $(t,x), (s,y)\in\mR^{1+nd}$,
\begin{align}\label{ES2}
\ell(t+s,x+y)\leq\ell(t,x)+\ell(s,y),
\end{align}
we have
$$
\ell(t-s,x-\e^{(t-s)A}y)\leq \ell(t-r,x-\e^{(t-r)A}z)+\ell(r-s,\e^{(t-r)A}z-\e^{(t-s)A}y).
$$
For simplicity, we write
$$
a:=\ell(t-r,x-\e^{(t-r)A}z),\ b:=\ell(r-s,z-\e^{(r-s)A}y).
$$
By the definition of $\ell$, we have
$$
|t-r|\leq a^2,\ \ |(z-\e^{(r-s)A}y)_j|\leq b^{1+2(n-j)},\ \ j=1,\cdots,n.
$$
Hence,  for each $i=1,\cdots,n$,
\begin{align}\label{EK3}
\begin{split}
&|(\e^{(t-r)A}z-\e^{(t-s)A}y)_i|=\left|\sum_{j=1}^n(\e^{(t-r)A})_{ij}(z-\e^{(r-s)A}y)_j\right|\\
&\quad\qquad\stackrel{\eqref{ES1}}{\leq}\sum_{j=i}^n\frac{|t-r|^{j-i}}{(j-i)!}b^{1+2(n-j)}\leq\sum_{j=i}^n\frac{a^{2(j-i)}}{(j-i)!}b^{1+2(n-j)}\\
&\quad\qquad\leq(a\vee b)^{1+2(n-i)}\sum_{j=i}^n\frac{1}{(j-i)!}\leq 3(a\vee b)^{1+2(n-i)},
\end{split}
\end{align}
and
$$
\ell(t-s,x-\e^{(t-s)A}y)\leq a+ 3(a\vee b)\leq 4(a+b).
$$
(iv) and (v) are easy consequences of (iii).
\end{proof}
For $f\in L^1_{loc}(\mR^{1+nd})$, we define the Hardy-Littlewood maximal function by
$$
\cM f(t,x):=\sup_{r>0}\fint_{Q_r(t,x)}|f(t',x')|\dif x'\dif t',
$$
and the sharp function by
$$
\cM^\sharp f(t,x):=\sup_{r>0}\fint_{Q_r(t,x)}|f(t',x')-f_{Q_r(t,x)}|\dif x'\dif t' .
$$
Here for a $Q\in\mQ$, 
$$
f_Q:=\fint_Q f(t',x')\dif x'\dif t' :=\frac{1}{|Q|}\int_Q f(t',x')\dif x'\dif t'.
$$
One says that a function $f\in BMO(\mR^{1+nd})$ if $\cM^\sharp f\in L^\infty(\mR^{1+nd})$. Clearly, $f\in BMO(\mR^{1+nd})$ if and only if
there exists a constant $C>0$ such that for any $Q\in\mQ$, and for some $c_Q\in\mR$,
$$
\fint_Q |f(t',x')-c_Q|\dif x'\dif t'\leq C.
$$
Using Lemma  \ref{Le31}, the following version of Fefferman-Stein type theorem can be established 
in a similar way as that for   \cite[Theorem 2.12]{Ch-Zh}.
We omit the details here.
\bt\label{Th2}
Suppose $q\in(1,\infty)$, and  $\sP$  is a
bounded linear operator from $L^q(\mR^{1+nd})$ to $L^q(\mR^{1+nd})$
and also from $L^\infty(\mR^{1+nd})$ to $BMO(\mR^{1+nd})$. 
Then for any $p\in[q,\infty)$, there is a  constant $C> 0$ depending only on $p,q$ and the norms of $\|\sP\|_{L^q\to L^q}$ and $\|\sP\|_{L^\infty\to BMO}$
so that 
$$
\|\sP f\|_p\leq C\|f\|_p  \qquad \hbox{for every } f\in L^p(\mR^{1+nd}) . 
$$
\et

\section{\bf Propagation of regularity in $L^2$-space for transport equations}

Fix $n\geq 2$ and $\lambda\geq 0$. Let $A$ be as in \eqref{AA}. 
In this section, corresponding to \eqref{e:1.12}, 
we consider the following linear transport equation in $\mR^{nd}$
for $u=u(s, x)$: 
\begin{align}\label{Tran}
\p_su+Ax\cdot\nabla u-\lambda u+f=0.
\end{align}
Taking Fourier's transform 
in the spatial variable 
 $x \in \mR^{nd}$, we obtain for any $(s, \xi)\in \mR \times \mR^{nd}$,  
\begin{align}\label{Fo}
\p_s\hat u (s, \xi) -A^*\xi\cdot\nabla_{\xi}  \hat u (x, \xi) -\lambda \hat u (s, \xi) +\hat f (s, \xi) =0,
\end{align}
where $A^*$ is the transpose of $A$. 
Multiplying both sides by the complex conjugate of $\hat u$, we get
$$
\p_s |\hat u|^2-A^*\xi\cdot\nabla|\hat u|^2-\lambda |\hat u|^2+2\Re(\hat f,\bar{\hat u})=0.
$$
Let $\phi(\xi)$ be a smooth function and define
\begin{align}\label{Gphi}
g_\phi:=2\phi\Re(\hat f,\bar{\hat u})+A^*\xi\cdot\nabla\phi|\hat u|^2.
\end{align}
We have
\begin{align}\label{E1}
\p_s (|\hat u|^2\phi)-A^*\xi\cdot\nabla(|\hat u|^2\phi)-\lambda (|\hat u|^2\phi)+g_\phi=0,
\end{align}
and if $\hat u$ has compact support,  then
\begin{align}\label{E2}
\begin{split}
(|\hat u|^2\phi)(s,\xi)&=-\int^{\infty}_0\p_t\Big(\e^{-\lambda t}(|\hat u|^2\phi)\big(t+s,\e^{-tA^*}\xi\big)\Big)\dif t\\
&=\int^{\infty}_0\e^{-\lambda t} g_\phi\big(t+s, \e^{-tA^*}\xi\big)\dif t,
\end{split}
\end{align}
where $\e^{-tA^*}=(\e^{-tA})^*$ is the transpose of exponential matrix in \eqref{ES1}.

\medskip

The following propagation of regularity in $L^2$-space 
is the key step in the proof of  Theorem \ref{Main},  which has independent interest.

\bt
Let $f, u\in L^2(\mR^{1+nd})$ so that \eqref{Tran} holds in the weak sense. 
For any $\alpha\geq 0$ and $j=1,2,\cdots,n-1$, there is a constant $C=C(\alpha,j,d)>0$ such that
\begin{align}\label{E4}
\|\Delta_{x_j}^{\frac{\alpha}{2(1+\alpha)}}u\|_2\leq C\|\Delta_{x_{j+1}}^{\frac{\alpha}{2}}u\|^{\frac{1}{1+\alpha}}_2\|f\|_2^{\frac{\alpha}{1+\alpha}}.
\end{align}
In particular,
\begin{align}\label{EV7}
\|\Delta_{x_j}^{\frac{\alpha}{2(1+(n-j)\alpha)}}u\|_2\leq C\|\Delta_{x_n}^{\frac{\alpha}{2}}u\|^{\frac{1+(n-j-1)\alpha}{1+(n-j)\alpha}}_2\|f\|_2^{\frac{\alpha}{1+(n-j)\alpha}}.
\end{align}
\et
\begin{proof}
Let $\rho: [0,\infty)\to[0,1]$ be a smooth function with $\rho (s)=1$ for $s<1$ and $\rho (s)=0$ for $s>2$. For $R>0$, let 
 $\eta_R(s,\xi):=\rho(|s|/R) \rho(|\xi|  /R  )$ 
and
$$
\hat u_R:=\eta_R\hat u, \ \ \hat f_R:=(A^*\xi\cdot\nabla\eta_R-\p_s\eta_R)\hat u+\eta_R\hat f.
$$
Since $\hat u_R$ satisfies
$$
\p_s\hat u_R-A^*\xi\cdot\nabla \hat u_R-\lambda\hat u_R+\hat f_R=0,
$$
if we can show that for some $C=C(\alpha,j,d)>0$,
$$
\|\Delta_{x_j}^{\frac{\alpha}{2(1+\alpha)}}u_R\|_2\leq C\|\Delta_{x_{j+1}}^{\frac{\alpha}{2}}u_R\|^{\frac{1}{1+\alpha}}_2\|f_R\|_2^{\frac{\alpha}{1+\alpha}},
$$ 
then letting $R\to\infty$, we get \eqref{E4}. Hence, without loss of generality, in the following, 
we may and do assume that $\hat u$ has compact support. We use   induction method.
\medskip\\
Let us first look at the case of $j=1$. We follow the simple argument of Alexander \cite{Al}. 
For any $\eps>0$, by Planchel's identity, we have
\begin{align*}
\|\Delta_{x_1}^{\frac{\alpha}{2(1+\alpha)}}u\|_2^2
&=\int^\infty_{-\infty}\int_{\mR^{nd}}|\xi_1|^{\frac{2\alpha}{1+\alpha}}|\hat u(s,\xi)|^2\dif \xi\dif s\\
&=\int^\infty_{-\infty}\int_{\mR^{nd}}|\xi_1|^{\frac{2\alpha}{1+\alpha}}1_{\{\eps|\xi_2|>|\xi_1|^{1/(1+\alpha)}\}}|\hat u(s,\xi)|^2\dif \xi\dif s\\
&+\int^\infty_{-\infty}\int_{\mR^{nd}}|\xi_1|^{\frac{2\alpha}{1+\alpha}}1_{\{\eps|\xi_2|\leq |\xi_1|^{1/(1+\alpha)}\}}|\hat u(s,\xi)|^2\dif \xi\dif s\\
&=:I_1(\eps)+I_2(\eps).
\end{align*}
For $I_1(\eps)$, we have
$$
I_1(\eps)\leq\eps^{2\alpha}\int^\infty_{-\infty}\int_{\mR^{nd}}|\xi_2|^{2\alpha}|\hat u(s,\xi)|^2\dif \xi\dif s=\eps^{2\alpha}\|\Delta_{x_2}^{\frac{\alpha}{2}}u\|^2_2.
$$
For $I_2(\eps)$, letting $\Omega:=\{\eps|\xi_2|\leq |\xi_1|^{1/(1+\alpha)}\}$ and by \eqref{E2} with $\phi=1$, we have
\begin{align*}
I_2(\eps)&=2\int^\infty_{-\infty}\int_{\mR^{nd}}|\xi_1|^{\frac{2\alpha}{1+\alpha}}1_\Omega(\xi)
\int^{\infty}_0\e^{-\lambda t} \Re(\hat f,\bar{\hat u})\big(t+s, \e^{-tA^*}\xi\big)\dif t\dif \xi\dif s\\
&=2\int^\infty_{-\infty}\int_{\mR^{nd}}\int^{\infty}_0 |(\e^{tA^*}\xi)_1|^{\frac{2\alpha}{1+\alpha}}1_\Omega(\e^{tA^*}\xi)\e^{-\lambda t}
\Re(\hat f,\bar{\hat u})\big(t+s, \xi\big)\dif t\dif \xi\dif s\\
&=2\int^\infty_{-\infty}\int_{\mR^{nd}}\left(\int^\infty_0|(\e^{tA^*}\xi)_1|^{\frac{2\alpha}{1+\alpha}}1_\Omega(\e^{tA^*}\xi)\e^{-\lambda t}\dif t\right)
\Re(\hat f,\bar{\hat u})\big(s, \xi\big)
 \dif \xi\dif s.
\end{align*}
Observing that
$$
\e^{tA^*}\xi=\left(\xi_1,t\xi_1+\xi_2,\tfrac{t^2}{2}\xi_1+t\xi_2+\xi_3,\cdots,\tfrac{t^{n-1}}{(n-1)!}\xi_{1}+\cdots+t\xi_{n-1}+\xi_n\right),
$$
since $\lambda\geq 0$, we have
\begin{align}\label{E3}
\begin{split}
&\int^{\infty}_0\!\!|(\e^{tA^*}\xi)_1|^{\frac{2\alpha}{1+\alpha}}1_{\Omega}(\e^{t A^*}\xi)\e^{-\lambda t}\dif t
=|\xi_1|^{\frac{2\alpha}{1+\alpha}}\int^{\infty}_0\!\!1_{\{\eps|\xi_2+t\xi_1|\leq |\xi_1|^{1/(1+\alpha)}\}}\e^{-\lambda t}\dif t\\
&\quad\leq|\xi_1|^{\frac{2\alpha}{1+\alpha}}\int^{\infty}_{-\infty}1_{\{|\xi_2|/|\xi_1|-|\xi_1|^{-\alpha/(1+\alpha)}/\eps
\leq t\leq |\xi_2|/|\xi_1|+|\xi_1|^{-\alpha/(1+\alpha)}/\eps\}}\dif t=2|\xi_1|^{\frac{\alpha}{1+\alpha}}/\eps.
\end{split}
\end{align}
Hence, by Young's inequality we have
\begin{align*}
I_2(\eps)&\leq\frac{4}{\eps}\int^\infty_{-\infty}\int_{\mR^{nd}}|\xi_1|^{\frac{\alpha}{1+\alpha}}(|\hat f|\,|\hat u|)(s, \xi)\dif\xi\dif s\\
&\leq\frac{1}{2}\int^\infty_{-\infty}\int_{\mR^{nd}}|\xi_1|^{\frac{2\alpha}{1+\alpha}}|\hat u|^2(s, \xi)\dif \xi\dif s
+\frac{8}{\eps^2}\int^\infty_{-\infty}\int_{\mR^{nd}}|\hat f|^2(s, \xi)\dif \xi\dif s\\
&=\frac{1}{2}\|\Delta_{x_1}^{\frac{\alpha}{2(1+\alpha)}}u\|^2_2+\frac{8}{\eps^2}\|f\|^2_2.
\end{align*}
Combining the above calculations, we obtain
$$
\|\Delta_{x_1}^{\frac{\alpha}{2(1+\alpha)}}u\|_2^2\leq
\eps^{2\alpha}\|\Delta_{x_2}^{\frac{\alpha}{2}}u\|^2_2+\frac{1}{2}\|\Delta_{x_1}^{\frac{\alpha}{2(1+\alpha)}}u\|^2_2+\frac{8}{\eps^2}\|f\|^2_2,
$$
which in turn gives \eqref{E4} for $j=1$ by letting $\eps=(\|f\|_2/\|\Delta^{\frac{\alpha}{2}}_{x_2} u\|_2)^{\frac{1}{1+\alpha}}$.
\medskip\\
Suppose now that \eqref{E4} has been proven for $j=1,\cdots,k$ with $k\leq n-2$. We want to show 
that \eqref{E4} holds for $j=k+1$.
For $\delta>0$, we define
$$
\phi_\delta(\xi):=\Pi_{j=1}^{k} \rho \Big(|\xi_j|  /(\delta^{k-j+1}|\xi_{k+1}|^{ (1+(k-j+2)\alpha)/(1+\alpha)})\Big),
$$
 and by Planchel's identity, write
\begin{align*}
\|\Delta_{x_{k+1}}^{\frac{\alpha}{2(1+\alpha)}}u\|_2
&=\int^\infty_{-\infty}\int_{\mR^{nd}}|\xi_{k+1}|^{\frac{2\alpha}{1+\alpha}}|\hat u(s,\xi)|^2\dif \xi\dif s\\
&=\int^\infty_{-\infty}\int_{\mR^{nd}}|\xi_{k+1}|^{\frac{2\alpha}{1+\alpha}}(1-\phi_\delta(\xi))|\hat u(s,\xi)|^2\dif \xi\dif s\\
&\quad+\int^\infty_{-\infty}\int_{\mR^{nd}}|\xi_{k+1}|^{\frac{2\alpha}{1+\alpha}}\phi_\delta(\xi)|\hat u(s,\xi)|^2\dif \xi\dif s\\
&=:K_1(\delta)+K_2(\delta).
\end{align*}
We first treat $K_1(\delta)$. By the induction hypothesis, one sees that for $j=1,2,\cdots,k$,
\begin{align}\label{E6}
\|\Delta_{x_j}^{\frac{\alpha}{2(1+(k-j+2)\alpha)}}u\|_2\leq 
C\|\Delta_{x_{k+1}}^{\frac{\alpha}{2(1+\alpha)}}u\|^{\frac{1+\alpha}{1+(k-j+2)\alpha}}_2\|f\|_2^{\frac{(k-j+1)\alpha}{1+(k-j+2)\alpha}}.
\end{align}
Observing that
$$
1-\phi_\delta(\xi)\leq\sum_{j=1}^{k}1_{\{|\xi_{j}|\geq\delta^{k-j+1}
 |\xi_{k+1}|^{(1+(k-j+2)\alpha)/(1+\alpha)}\}},
$$
by \eqref{E6} and Young's inequality, we have
\begin{align*}
K_1(\delta)&\leq\sum_{j=1}^{k}
\int^\infty_{-\infty}\int_{\mR^{nd}}|\xi_{k+1}|^{\frac{2\alpha}{1+\alpha}}1_{\{|\xi_{j}|\geq\delta^{k-j+1}
 |\xi_{k+1}|^{(1+(k-j+2)\alpha)/(1+\alpha)}\}}|\hat u(s,\xi)|^2\dif \xi\dif s \\
&\leq\sum_{j=1}^{k}\delta^{-\frac{2(k-j+1)\alpha}{1+(k-j+2)\alpha}}\int^\infty_{-\infty}\int_{\mR^{nd}}
|\xi_j|^{\frac{2\alpha}{1+(k-j+2)\alpha}}|\hat u(s,\xi)|^2\dif \xi\dif s
=\sum_{j=1}^{k}\delta^{-\frac{2(k-j+1)\alpha}{1+(k-j+2)\alpha}}\|\Delta^{\frac{\alpha}{2(1+(k-j+2)\alpha)}}_{x_j}u\|^2_2\\
 &\leq C\sum_{j=1}^{k}\|\Delta_{x_k}^{\frac{\alpha}{2(1+\alpha)}}u\|^{\frac{2(1+\alpha)}{1+(k-j+2)\alpha}}_2
(\delta^{-1}\|f\|_2)^{\frac{2(k-j+1)\alpha}{1+(k-j+2)\alpha}}
\leq\frac{1}{4}\|\Delta_{x_{k+1}}^{\frac{\alpha}{2(1+\alpha)}}u\|^2_2+C\delta^{-2}\|f\|_2^2.
\end{align*}
Next we treat the trouble term $K_2(\delta)$. Let 
$$
\xi_{(k)}:=(\xi_1,\xi_{2},\cdots,\xi_{k}), \ \xi^{(k)}:=(\xi_{k+1},\xi_{k+2},\cdots,\xi_n)
$$ 
and define
$$
h_\delta(s,\xi^{(k)}):=\int_{\mR^{kd}}(|\hat u|^2\phi_\delta)(s,\xi_{(k)},\xi^{(k)})\dif\xi_{(k)}.
$$
Integrating both sides of \eqref{E1} 
with respect to the first $k$-variables $\xi_{(k)}$,
we get
$$
\p_s h_\delta-A_{n-k}\xi^{(k)}\cdot\nabla_{\xi^{(k)}} h_\delta-\lambda h_\delta
=\int_{\mR^{kd}}\left(\sum_{j=1}^{k}\xi_j\cdot\nabla_{\xi_{j+1}}(|\hat u|^2\phi_{\delta})-g_{\phi_\delta}\right)\dif\xi_{(k)}=:g^{(k)}_{\phi_\delta},
$$
 and by \eqref{E2},
\begin{align*}
h_\delta(s,\xi^{(k)})=\int^{\infty}_0\e^{-\lambda t}g^{(k)}_{\phi_\delta}(t+s,\e^{-tA^*_{n-k}}\xi^{(k)})\dif t.
\end{align*}
Let
$$
\psi_\eps(\xi^{(k)}):=\rho \Big(\eps|\xi^{(k)}_{2}|/|\xi^{(k)}_{1}|^{1/(1+\alpha)}\Big)=\rho \Big(\eps|\xi_{k+2}|/|\xi_{k+1}|^{1/(1+\alpha)}\Big).
$$
We make the following decomposition as above:
\begin{align*}
K_2(\delta)&=\int^\infty_{-\infty}\!\int_{\mR^{(n-k)d}}|\xi_{k+1}|^{\frac{2\alpha}{1+\alpha}}\left(1-\psi_\eps(\xi^{(k)})\right)h_\delta(s,\xi^{(k)})\dif\xi^{(k)}\dif s\\
&\quad+\int^\infty_{-\infty}\!\int_{\mR^{(n-k)d}}|\xi_{k+1}|^{\frac{2\alpha}{1+\alpha}}\psi_\eps(\xi^{(k)})h_\delta(s,\xi^{(k)})\dif\xi^{(k)}\dif s\\
&=:K_{21}(\delta,\eps)+K_{22}(\delta,\eps).
\end{align*}
For $K_{21}(\delta,\eps)$, we have
\begin{align*}
K_{21}(\delta,\eps)&\leq\int^\infty_{-\infty}\!\int_{\mR^{(n-k)d}}|\xi_{k+1}|^{\frac{2\alpha}{1+\alpha}}
1_{\{\eps|\xi_{k+2}|>|\xi_{k+1}|^{1/(1+\alpha)}\}}h_\delta(s,\xi^{(k)})\dif \xi^{(k)}\dif s\\
&\leq\eps^{2\alpha}\int^\infty_{-\infty}\!\int_{\mR^{nd}}|\xi_{k+2}|^{2\alpha}|\hat u|^2(s,\xi)\dif\xi\dif s
=\eps^{2\alpha}\|\Delta^{\frac{\alpha}{2}}_{x_{k+2}}u\|^2_2.
\end{align*}
For $K_{22}(\delta,\eps)$, 
 by the change of variables and Fubini's theorem,
we have
\begin{align*}
&K_{22}(\delta,\eps)=\int^\infty_{-\infty}\!\int_{\mR^{(n-k)d}}|\xi^{(k)}_{1}|^{\frac{2\alpha}{1+\alpha}}\psi_\eps(\xi^{(k)})
\int^{\infty}_0\e^{-\lambda t}g^{(k)}_\phi(t+s,\e^{-tA^*_{n-k}}\xi^{(k)})\dif t\dif \xi^{(k)}\dif s\\
&\quad=\int^\infty_{-\infty}\int_{\mR^{(n-k)d}}\left(\int^\infty_0|(\e^{tA^*_{n-k}}\xi^{(k)})_{1}|^{\frac{2\alpha}{1+\alpha}}
\psi_\eps(\e^{tA^*_{n-k}}\xi^{(k)})\e^{-\lambda t}\dif t\right)g^{(k)}_\phi(s,\xi^{(k)})\dif \xi^{(k)}\dif s.
\end{align*}
Letting
\begin{align*}
\gamma_\eps(\xi^{(k)})&:=\int^\infty_0|(\e^{tA^*_{n-k}}\xi^{(k)})_{1}|^{\frac{2\alpha}{1+\alpha}}\psi_\eps(\e^{tA^*_{n-k}}\xi^{(k)})\e^{-\lambda t}\dif t\\
&=|\xi_{k+1}|^{\frac{2\alpha}{1+\alpha}}\int^\infty_0
\rho (\eps|\xi_{k+2}+t\xi_{k+1}|/|\xi_{k+1}|^{1/(1+\alpha)})\e^{-\lambda t}\dif t,
\end{align*}
and recalling
$$
g^{(k)}_{\phi_\delta}=
 -\int_{\mR^{kd}}g_{\phi_\delta}(\xi_{(k)},\cdot)\dif\xi_{(k)}
+\sum_{j=1}^k\int_{\mR^{kd}}\xi_j\cdot\nabla_{\xi_{j+1}}(|\hat u|^2\phi_\delta)\dif\xi_{(k)},
$$
we may write
\begin{align*}
K_{22}(\delta,\eps)&=
 -\int^\infty_{-\infty}\int_{\mR^{(n-k)d}}\gamma_\eps(\xi^{(k)})\int_{\mR^d}g_{\phi_\delta}(s,\xi_{(k)},\xi^{(k)})\dif\xi_{(k)}\dif \xi^{(k)}\dif s\\
&+\sum_{j=1}^k\int^\infty_{-\infty}\int_{\mR^{(n-k)d}}\gamma_\eps(\xi^{(k)})\int_{\mR^d}\xi_j\cdot\nabla_{\xi_{j+1}}(|\hat u|^2\phi_\delta)\dif\xi_{(k)}\dif \xi^{(k)}\dif s\\
&=:K_{221}(\delta,\eps)+K_{222}(\delta,\eps).
\end{align*}
As in estimating \eqref{E3}, we have
$$
\gamma_\eps(\xi^{(k)})\leq 4|\xi_{k+1}|^{\frac{\alpha}{1+\alpha}}/\eps,
$$
and also by the definition of $g_{\phi_\delta}$,
$$
 |g_{\phi_\delta}|\leq 2|\hat f|\,|\hat u|+\sum_{j=1}^{n-1}|\xi_j|\cdot|\nabla_{j+1}\phi_\delta|\cdot|\hat u|^2
\leq 2|\hat f|\,|\hat u|+C\delta|\xi_{k+1}|^{\frac{\alpha}{1+\alpha}}|\hat u|^2.
$$
Hence,
\begin{align*}
K_{221}(\delta,\eps)&\leq\frac{8}{\eps}\int^\infty_{-\infty}\int_{\mR^{nd}}|\xi_{k+1}|^{\frac{\alpha}{1+\alpha}}|\hat f|\,|\hat u|+
\frac{C\delta}{\eps}\int^\infty_{-\infty}\int_{\mR^{nd}}|\xi_{k+1}|^{\frac{2\alpha}{1+\alpha}}|\hat u|^2\\
&\leq\left(\tfrac{1}{4}+\tfrac{C\delta}{\eps}\right)\|\Delta_{x_{k+1}}^{\frac{\alpha}{2(1+\alpha)}}u\|^2_2+C\eps^{-2}\|f\|_2^2.
\end{align*}
 On the other hand, by elementary calculations, 
we also have 
\begin{align*}
|\nabla_{\xi_{k+1}}\gamma_\eps(\xi^{(k)})|\leq C|\xi_{k+1}|^{-\frac{1}{1+\alpha}}/\eps.
\end{align*}
Thus, since $\hat u$ has compact support, by the integration by parts formula, we have
\begin{align*}
|K_{222}(\delta,\eps)|&=\left|\int^\infty_{-\infty}\int_{\mR^{(n-k)d}}\!\int_{\mR^{kd}}
 \xi_k\cdot\nabla_{\xi_{k+1}}\gamma_\eps(\xi^{(k)})\,(|\hat u|^2\phi_\delta)\dif\xi_{(k)}\dif \xi^{(k)}\dif s\right|\\
&\leq\frac{C}{\eps}\int^\infty_{-\infty}\int_{\mR^{nd}}|\xi_k|\,|\xi_{k+1}|^{-\frac{1}{1+\alpha}}(|\hat u|^2\phi_\delta)\dif \xi\dif s\\
&\leq\frac{C}{\eps}\int^\infty_{-\infty}\int_{\mR^{nd}}|\xi_k|\,|\xi_{k+1}|^{-\frac{1}{1+\alpha}}|\hat u|^21_{\{|\xi_{k}|\leq\delta|\xi_{k+1}|^{(1+2\alpha)/(1+\alpha)}\}}\dif \xi\dif s\\
&\leq\frac{C\delta}{\eps}\int^\infty_{-\infty}\int_{\mR^{nd}}|\xi_{k+1}|^{\frac{2\alpha}{1+\alpha}}|\hat u|^2\dif s
=\frac{C\delta}{\eps}\|\Delta_{x_{k+1}}^{\frac{\alpha}{2(1+\alpha)}}u\|^2_2.
\end{align*}
Combining the above calculations, we obtain
$$
\|\Delta_{x_{k+1}}^{\frac{\alpha}{2(1+\alpha)}}u\|^2_2\leq \left(\tfrac{1}{4}+\tfrac{C_1'\delta}{\eps}\right)\|\Delta_{x_{k+1}}^{\frac{\alpha}{2(1+\alpha)}}u\|^2_2
+\eps^{2\alpha}\|\Delta^{\frac{\alpha}{2}}_{x_{k+2}}u\|^2_2+C_2'(\eps^{-2}+\delta^{-2})\|f\|_2^2.
$$
Choosing $\delta=\eps/(4C_1')$, we get
$$
\|\Delta_{x_{k+1}}^{\frac{\alpha}{2(1+\alpha)}}u\|^2_2\leq 2\eps^{2\alpha}\|\Delta^{\frac{\alpha}{2}}_{x_{k+2}}u\|^2_2+2C'_2(1+16(C'_1)^2)\eps^{-2}\|f\|_2^2,
$$
which then yields \eqref{E4} for $j=k+1$  by letting $\eps=(\|f\|_2/\|\Delta^{\frac{\alpha}{2}}_{x_{k+2}} u\|_2)^{\frac{1}{1+\alpha}}$.
\end{proof}

\section{\bf Proof of Theorem \ref{Main}}

\subsection{Case: $p=2$}
Without loss of generality, we may assume $f\in C^\infty_c(\mR^{1+nd})$.
It follows from  Fourier's transform
and H\"older's inequality that
\begin{align*}
 &\int^\infty_{-\infty}\left\|\Delta_{x_n}\int^\infty_0\e^{-\lambda t}\cT_{s,t+s}f(t+s,\cdot)\dif t\right\|_2^2\dif s\\
&=\int^\infty_{-\infty}\int_{\mR^{nd}}|\xi_n|^{2}\left|\int^\infty_0\e^{-\lambda t}\widehat{\cT_{s,t+s}f}(t+s,\xi)\dif t\right|^2\dif\xi\dif s\\
&=\int^\infty_{-\infty}\int_{\mR^{nd}}|\xi_n|^{2}\left|\int^\infty_0\e^{-\lambda t}
\e^{-\frac{1}{2}\int^t_0|(\sigma^a_{s+r})^*\e^{-rA^*}\xi|^2\dif r}\hat f(t+s,\e^{-tA^*}\xi )\dif t\right|^2\dif\xi\dif s\\
&\leq\int^\infty_{-\infty}\int_{\mR^{nd}}\left(\int^\infty_0|\xi_n|^2
\e^{-\frac{1}{2}\int^t_0|(\sigma^a_{s+r})^*\e^{-rA^*}\xi|^2\dif r}|\hat f(t+s,\e^{-tA^*}\xi)|^2\dif t\right)\\
&\qquad\qquad\times\left(\int^\infty_0|\xi_n|^2\e^{-\frac{1}{2}\int^t_0|(\sigma^a_{s+r})^*\e^{-rA^*}\xi|^2\dif r}\dif t\right)\dif\xi\dif s.
\end{align*}
Let $\Theta_{t}$ be the dilation operator defined by \eqref{Dia}. By definitions \eqref{ES11}, \eqref{ES1} and \eqref{Dia}, it is easy to see that
$$
t^{1/2}(\sigma^a_u)^*\e^{-rtA^*}\xi=(\sigma^a_u)^*\e^{-rA^*}\Theta_{t^{1/2}}(\xi).
$$
Thus, by the change of variables and as in showing \eqref{HJ2}, we have for some $c>0$,
$$
\int^t_0|(\sigma^a_{s+r})^*\e^{-rA^*}\xi|^2\dif r=\int^1_0|(\sigma^a_{s+tr})^*\e^{-rA^*}\Theta_{t^{1/2}}(\xi)|^2\dif r\geq c|\Theta_{t^{1/2}}\xi|^2,
$$
and similarly,
$$
\int^t_0|(\sigma^a_{s+r})^*\e^{(t-r)A^*}\xi|^2\dif r\geq c|\Theta_{t^{1/2}}\xi|^2.
$$
Hence,
\begin{align*}
\int^\infty_0|\xi_n|^2&\e^{-\frac{1}{2}\int^t_0|(\sigma^a_{s+r})^*\e^{-rA^*}\xi|^2\dif r}\dif t
\leq\int^\infty_0|\xi_n|^2\e^{-c|\Theta_{t^{1/2}}\xi|^2}\dif t\leq\int^\infty_0|\xi_n|^2\e^{-ct|\xi_n|^2}\dif t=c^{-1},
\end{align*}
and
\begin{align*}
&\int^{\infty}_0|(\e^{tA^*}\xi)_n|^2
\e^{-\frac{1}{2}\int^t_0|(\sigma^a_{s-t+r})^*\e^{(t-r)A^*}\xi|^2\dif r}\dif t
\leq\int^{\infty}_0|(\e^{tA^*}\xi)_n|^2\e^{-c|\Theta_{t^{1/2}}\xi|^2}\dif t\\
&\qquad\leq\sum_{j=1}^n\frac{1}{(n-j)!}\int^\infty_0|t^{n-j}\xi_j|^2\e^{-ct^{2(n-j)+1}|\xi_j|^2}\dif t
\leq 2c^{-1}.
\end{align*}
Thus, by the change of variables and Fubini's theorem, we further have
\begin{align*}
&\int^\infty_{-\infty}\left\|\Delta_{x_n}\int^\infty_0\e^{-\lambda t}\cT_{s,t+s}f(t+s,\cdot)\dif t\right\|_2^2\dif s\\
&\leq c^{-1}\int^\infty_{-\infty}\int_{\mR^{nd}}\left(\int^\infty_0|\xi_n|^2
\e^{-\frac{1}{2}\int^t_0|(\sigma^a_{s+r})^*\e^{-rA^*}\xi|^2\dif r}|\hat f(t+s,\e^{-tA^*}\xi)|^2\dif t\right)\dif\xi\dif s\\
&=c^{-1}\int^\infty_{-\infty}\int_{\mR^{nd}}\left(\int^\infty_0|(\e^{tA^*}\xi)_n|^2
\e^{-\frac{1}{2}\int^t_0|(\sigma^a_{s+r})^*\e^{(t-r)A^*}\xi|^2\dif r}|\hat f(t+s,\xi)|^2\dif t\right)\dif\xi\dif s\\
&=c^{-1}\int^\infty_{-\infty}\int_{\mR^{nd}}\left(\int^{\infty}_0|(\e^{tA^*}\xi)_n|^2
\e^{-\frac{1}{2}\int^t_0|(\sigma^a_{s-t+r})^*\e^{(t-r)A^*}\xi|^2\dif r}\dif t\right)|\hat f(s,\xi)|^2\dif\xi\dif s\\
&\leq 2c^{-2}\int^\infty_{-\infty}\!\int_{\mR^{nd}}|\hat f(s,\xi)|^2\dif\xi\dif s=2c^{-2}\|f\|^2_2.
\end{align*}
This together with \eqref{EV7} completes the proof of \eqref{BN5} for $p=2$.

\subsection{Case: $p\in(2,\infty)$}
 Let $\varrho\in C^\infty_c(\mR^{nd})$ be nonnegative  with   $\int\varrho=1$.
We use it to define a family of
 mollifiers 
  $$
 \varrho_\eps(x)=\eps^{-nd}\varrho(x/\eps);\ \ \eps >0.
 $$ 
 For a function $f(t, x)$ defined on $\mR \times \mR^{nd}$ and $\eps >0$, let
$$ 
f_\eps(t,x) := f (t, \cdot )* \varrho_\eps (x):=\int_{\mR^d} f(t, y) \varrho_\eps(x-y ) \dif y .
$$
For $j=1,\cdots,n$ and $\eps\in(0,1)$,   define
\begin{align*}
\sP^\eps_jf:=\sP^{a}_j f_\eps(s,x)&:=\Delta_{x_j}^{ {1}/{(1+2(n-j))}}\int^\infty_0\e^{-\lambda t}\cT^{a}_{s,t+s} f_\eps(t+s,x)\dif t,
\end{align*}
where   the superscript $a$ denotes the dependence on the diffusion coefficient $a$.
To use Theorem \ref{Th2}, our main task is to show that $\sP^\eps_j$ is a bounded linear operator from $L^\infty(\mR^{1+nd})$ to $BMO$.
More precisely,  we want to prove that for any $f\in L^\infty(\mR^{1+nd})$ with $\|f\|_\infty\leq 1$, and any $Q=Q_r(t_0,x_0)\in\mQ$,
\begin{align}
\fint_{Q}|\sP^a_j f_\eps(s,x)-c^{Q}_j|^2 \leq C,\label{ET0}
\end{align}
where $c^{Q}_j$ is a constant depending on $Q$ and $f_\eps$, and $C$ only depends on $n,\kappa, p,d$, and not on $\eps$.
\bl\label{Le34}
(Scaling Property) For any $Q=Q_r(t_0,x_0)\in\mQ$ , we have
\begin{align}
\fint_{Q_r(t_0,x_0)}\big|\sP^a_jf_\eps(s,x)-c\big|^2=
\fint_{Q_1(0)}\big|\sP^{\tilde a}_j\tilde f_\eps(s,x)-c\big|^2,\label{SCL}
\end{align}
where $c\in\mR$, $\tilde a_s:=a_{r^2 s+t_0}$ and $\tilde f_\eps(t,x):=f_\eps\big(r^2 t+t_0, \Theta_r x+\e^{tA}x_0\big).$ Here 
$\Theta_r$ is the dilation operator defined in \eqref{Dia}.
\el
\begin{proof}
Let us write
$$
u_\eps(s,x):=\int^\infty_0\e^{-\lambda t}\cT^{a}_{s,t+s} f_\eps(t+s,x)\dif t.
$$
By the change of variables, we have
$$
\tilde u_\eps(s,x):=r^{-2} u_\eps\big(r^2 s+t_0,\Theta_rx+\e^{s A}x_0\big)=\int^\infty_0\e^{-\lambda t}\cT^{\tilde a}_{s,t+s} \tilde f_\eps(t+s,x)\dif t.
$$
Noticing that
$$
(\Delta_{x_j}^{ {1}/{(1+2(n-j))}}\tilde u_\eps)(s,x)=(\Delta_{x_j}^{ {1}/{(1+2(n-j))}} u_\eps)\big(r^2 s+t_0,\Theta_rx+\e^{s A}x_0\big),
$$
by the change of variables again, we have
$$
\fint_{Q_r(t_0,x_0)}\big|  \Delta_{x_j}^{ {1}/{(1+2(n-j))}} u_\eps(s,x)-c\big|^2
=\fint_{Q_1(0)}\big| \Delta_{x_j}^{ {1}/{(1+2(n-j))}} \tilde u_\eps(s,x)-c\big|^2.
$$
The proof is finished.
\end{proof}

Noticing that
$$
\int^\infty_0\e^{-\lambda t}\cT^{a}_{s,t+s} f_\eps(t+s,x)\dif t=\int^\infty_s\e^{\lambda (s-t)}\cT^{a}_{s,t} f_\eps(t,x)\dif t,
$$
for $s\in[-1,1]$, we split $\sP^\eps_jf(s,x)=\sP^\eps_{j1}f(s,x)+\sP^\eps_{j2}f(s,x)$ with
\begin{align*}
\sP^\eps_{j1} f(s,x)&:= \Delta_{x_j}^{ {1}/{(1+2(n-j))}} \int^2_s\e^{\lambda(s-t)}\cT_{s,t}f_\eps(t,x)\dif t, \\
\sP^\eps_{j2} f(s,x)&:= \Delta_{x_j}^{ {1}/{(1+2(n-j))}} \int^\infty_2\e^{\lambda(s-t)}\cT_{s,t} f_\eps(t,x)\dif t.
\end{align*}
In the rest of this paper,   unless otherwise specified, 
all the constants contained in ``$\lesssim$'' will depend only on $n,\kappa,p,d$.

\bl\label{Le35}
Under  \eqref{Ell}, there is a constant $C=C(n,\kappa,p,d)>0$
such that for all $f\in L^\infty(\mR^{1+nd})$ with $\|f\|_\infty\leq 1$,
\begin{align}
\sup_{\eps\in(0,1)}\int_{Q_1(0)}|\sP^\eps_{j1} f(s,x)|^2\leq C.\label{BN4}
\end{align}
\el
\begin{proof}
For $s\in[-1,1]$, let
$$
u_\eps(s,x):=\int^2_s\e^{\lambda(s-t)}\cT_{s,t}f_\eps (t,x)\dif t=\int^\infty_s\e^{\lambda(s-t)}\cT_{s,t}((1_{[-1,2]}f_\eps) (t))(x)\dif t.
$$
Since $\|f\|_\infty\leq 1$, we have
\begin{align}
\|u_\eps(s)\|_\infty\leq 3,\ \ s\in[-1,1].\label{BOU}
\end{align}
By \eqref{EK2},  we have for any $s\in[-1,1]$,
\begin{align}\label{EK5}
\|\nabla_{x_n} u_\eps(s)\|_\infty\leq\int^2_s\|\nabla_{x_n}\cT_{s,t}f_\eps (t)\|_\infty\dif t
\lesssim\int^2_s(t-s)^{-1/2}\dif t\lesssim 1.
\end{align}
Let $\varphi$ be a nonnegative smooth cutoff function in $\mR^{nd}$ with $\varphi(x)=1$ for $|x|\leq n$ and  $\varphi(x)=0$ for $|x|>2n$.
Notice that
$$
\p_s u_\eps+\sL_s u_\eps+f_\eps1_{[-1,2]}(s)=0.
$$
It is easy to see that $u_\eps\varphi$ satisfies
$$
\p_s(u_\eps\varphi)+\sL_s (u_\eps\varphi)=\sL_s(u_\eps\varphi)-(\sL_s u_\eps)\varphi-f_\eps\varphi1_{[-1,2]}(s)=:g^\varphi_\eps,
$$
which implies by \eqref{Gen} that
$$
(u_\eps\varphi)(s,x)=\int^\infty_s\e^{\lambda(s-t)}\cT_{s,t}g^\varphi_\eps(t,x)\dif t.
$$
By the definition of $\sP^\eps_{j1}$, we have
\begin{align*}
&\int_{Q_1(0)}|\sP^\eps_{j1} f|^2=\int_{Q_1(0)}| \Delta_{x_j}^{ {1}/{(1+2(n-j))}}u_\eps|^2\leq\int_{\mR^{1+nd}}| (\Delta_{x_j}^{ {1}/{(1+2(n-j))}} u_\eps)\varphi|^2\\
&\quad\leq 2\big\| \Delta_{x_j}^{ {1}/{(1+2(n-j))}}(u_\eps\varphi)\big\|_2^2+
2\big\| \Delta_{x_j}^{ {1}/{(1+2(n-j))}} (u_\eps\varphi)- (\Delta_{x_j}^{ {1}/{(1+2(n-j))}} u_\eps)\varphi\big\|_2^2=:I_1+I_2.
\end{align*}
For $I_1$,  by \eqref{BN5} for $p=2$ and \eqref{BOU}, \eqref{EK5}, we have
$$
I_1=\left\| \Delta_{x_j}^{ {1}/{(1+2(n-j))}}\int^\infty_s\e^{\lambda(s-t)}\cT_{s,t}g^\varphi_\eps(t,x)\dif t\right\|_2\lesssim \|g^\varphi_\eps\|_2\lesssim 1.
$$
For $I_2$, if $j=n$, then by \eqref{BOU} and \eqref{EK5},
$$
I_2=2\big\|u_\eps\Delta_{x_n}\varphi+2\nabla_{x_n}u_\eps\cdot\nabla_{x_n}\varphi\big\|_2^2\lesssim 1;
$$
if $j=1,\cdots,n-1$, then by definition \eqref{Loc} and \eqref{BOU},
\begin{align*}
I_2&=\left\|\int_{\mR^d}\frac{\delta^{(1)}_{z_j} u_\eps\,\delta^{(1)}_{z_j} \varphi}{|z_j|^{d+2/(1+2(n-j))}}\dif z_j\right\|_2
\leq \int_{\mR^d}\frac{\|\delta^{(1)}_{z_j} u_\eps\,\delta^{(1)}_{z_j} \varphi\|_2}{|z_j|^{d+2/(1+2(n-j))}}\dif z_j\\
&\leq 6\int_{\mR^d}\frac{\|\delta^{(1)}_{z_j} \varphi\|_2}{|z_j|^{d+2/(1+2(n-j))}}\dif z_j\lesssim\int_{\mR^d}\frac{1\wedge|z_j|}{|z_j|^{d+2/(1+2(n-j))}}\dif z_j\lesssim 1.
\end{align*}
The proof is complete.
\end{proof}

The next lemma is crucial for treating $\sP^\eps_{j2}f$.
\bl
Under  \eqref{AA}, there is a constant $C=C(n,\kappa, p,d)>0$
such that for all $f\in L^\infty(\mR^{1+nd})$ with $\|f\|_\infty\leq 1$ and all $s\in[-1,1]$,
\begin{align}
\sup_{\eps\in(0,1)}\int^\infty_2\Big| \Delta_{x_j}^{ {1}/{(1+2(n-j))}} \cT_{s,t}f_\eps(t,0)
-  \Delta_{x_j}^{ {1}/{(1+2(n-j))}} \cT_{0,t}f_\eps(t,0)\Big|\dif t\leq C.\label{NB6}
\end{align}
\el
\begin{proof}
Noticing that by \eqref{Gen}, for each $s\in(0,t)$ and $x\in\mR^{nd}$,
\begin{align}\label{In}
\cT_{0,t}f_\eps(t,x)-\cT_{s,t}f_\eps(t,x)=\int^s_0\sL_r\cT_{r,t}f_\eps(t,x)\dif r,
\end{align}
we have
\begin{align*}
I^\eps_j(s,t)&:=\Delta_{x_j}^{ {1}/{(1+2(n-j))}}\cT_{0,t}f_\eps(t,0)
-\Delta_{x_j}^{ {1}/{(1+2(n-j))}}\cT_{s,t}f_\eps(t,0)\\
&=\int^s_0(\Delta_{x_j}^{ {1}/{(1+2(n-j))}}\sL_r\cT_{r,t}f_\eps)(t,0)\dif r.
\end{align*}
For $j=1$, we have for all $s\in[-1,1]$,
\begin{align*}
\int^\infty_2|I^\eps_1(s,t)|\dif t
&=\int^\infty_2\left|\int^s_0\Big(\Delta_{x_1}^{\frac{1}{1+2(n-1)}}\tr(a_r\cdot\nabla^2_{x_n})\cT_{r,t}f_\eps\Big)(t,0)\dif r\right|\dif t\\
&\leq\kappa\int^\infty_2\!\!\!\int^s_0|\Delta_{x_1}^{\frac{1}{1+2(n-1)}}\nabla^2_{x_n}\cT_{r,t}f_\eps|(t,0)\dif r\dif t\\
&\stackrel{\eqref{EK2}}{\lesssim} \int^\infty_2\!\!\int^s_0(t-r)^{-3/2}\dif r\dif t\lesssim 1.
\end{align*}
For $j=n$, we have for all $s\in[-1,1]$,
\begin{align*}
\int^\infty_2\!\!\!\!|I^\eps_n(s,t)|\dif t
&=\int^\infty_2\left|\int^s_0\Big(\big(\Delta_{x_n}\tr(a_r\cdot\nabla^2_{x_n})+(\nabla_{x_n}\cdot\nabla_{x_{n-1}})\big)\cT_{r,t}f_\eps\Big)(t,0)\dif r\right|\dif t\\
&\leq\int^\infty_2\!\!\!\int^s_0\Big(\kappa|\Delta_{x_n}\nabla^2_{x_n}\cT_{r,t}f_\eps|
+|(\nabla_{x_n}\cdot\nabla_{x_{n-1}})\cT_{r,t}f_\eps|\Big)(t,0)\dif r\dif t\\
&\stackrel{\eqref{EK2}}{\lesssim}\int^\infty_2\!\!\int^s_0(t-r)^{-2}\dif r\dif t\lesssim 1.
\end{align*}
For $j=2,\cdots, n-1$, since $\Delta_{x_j}^{ {1}/{(1+2(n-j))}}$ is a nonlocal operator, 
we have to carefully treat the trouble term $\Delta_{x_j}^{ {1}/{(1+2(n-j))}}(x_j\cdot\nabla_{x_{j-1}}\cT_{r,t}f_\eps)(t,0)$.
Fix
$$
\gamma\in\Big(\tfrac{2(n-j)+1}{2},\tfrac{(2(n-j)+1)^2}{4(n-j)}\Big).
$$
By \eqref{Loc} and \eqref{In}, we may write
\begin{align*}
I^\eps_j(s,t)&=\int_{|z_j|>t^\gamma}\Big(\delta_{z_j}\cT_{0,t}f_\eps(t,0)-\delta_{z_j}\cT_{s,t}f_\eps(t,0)\Big)\frac{\dif z_j}{|z_j|^{d+2/(1+2(n-j))}}\\
&\quad+\int_{|z_j|\leq t^\gamma}\left(\int^s_0\delta_{z_j}\tr(a_r\cdot\nabla^2_{x_n})\cT_{r,t}f_\eps(t,0)\dif r\right)\frac{\dif z_j}{|z_j|^{d+2/(1+2(n-j))}}\\
&\quad+\int_{|z_j|\leq t^\gamma}\left(\int^s_0\delta_{z_j}(Ax\cdot\nabla\cT_{r,t}f_\eps)(t,0)\dif r\right)\frac{\dif z_j}{|z_j|^{d+2/(1+2(n-j))}}\\
&=:I^\eps_{j1}(s,t)+I^\eps_{j2}(s,t)+I^\eps_{j3}(s,t).
\end{align*}
For $I^\eps_{j1}(s,t)$, thanks to $\gamma>\tfrac{2(n-j)+1}{2}$, we have for all $s\in[-1,1]$,
\begin{align*}
\int^\infty_2|I^\eps_{j1}(s,t)|\dif t
\lesssim \int^\infty_2\!\!\!\int_{|z_j|>t^\gamma}\frac{\dif z_j}{|z_j|^{d+2/(1+2(n-j))}}\dif t\lesssim\int^\infty_2|t|^{-2\gamma/(1+2(n-j))}\dif t\lesssim 1.
\end{align*}
For $I^\eps_{j2}(s,t)$, by \eqref{EK2} and $\gamma<\tfrac{(2(n-j)+1)^2}{2(2(n-j)-1)}$,
we have for all $s\in[-1,1]$,
\begin{align*}
\int^\infty_2|I^\eps_{j2}(s,t)|\dif t&\lesssim \int^\infty_2\!\!\!\int^s_0
\left(\int_{|z_j|\leq t^\gamma}|z_j|\,\|\nabla_{x_j}\nabla^2_{x_n}\cT_{r,t}f_\eps\|_\infty\frac{\dif z_j}{|z_j|^{d+2/(1+2(n-j))}}\right)\dif r\dif t\\
&\lesssim \int^\infty_2\!\!\int^s_0(t-r)^{-(2(n-j)+3)/2}\left(\int_{|z_j|\leq t^\gamma}\frac{|z_j|\dif z_j}{|z_j|^{d+2/(1+2(n-j))}}\right)\dif r\dif t\\
&\lesssim\int^\infty_2(t-1)^{-(2(n-j)+3)/2} t^{\gamma(1-2/(1+2(n-j)))}\dif t\lesssim 1.
\end{align*}
For $I^\eps_{j3}(s,t)$,  letting $\tilde z_j=(0,\cdots, 0,z_j,0,\cdots,0)$ and observing that
$$
\int_{|z_j|\leq t^\gamma}\frac{z_j\dif z_j}{|z_j|^{d+2/(1+2(n-j))}}=0,
$$
by \eqref{EK2} and $\gamma<\tfrac{(2(n-j)+1)^2}{4(n-j)}$, we have for all $s\in[-1,1]$,
\begin{align*}
&\int^\infty_2\!\!\!|I^\eps_{j3}(s,t)|\dif t=\int^\infty_2\left|\int^s_0\!\!\!\int_{|z_j|\leq t^\gamma}z_j\cdot\nabla_{x_{j-1}}
\cT_{r,t}f_\eps(t,\tilde  z_j)\frac{\dif z_j}{|z_j|^{d+2/(1+2(n-j))}}\dif r\right|\dif t\\
&=\int^\infty_2\left|\int^s_0\!\!\!\int_{|z_j|\leq t^\gamma}z_j\cdot\nabla_{x_{j-1}}(\cT_{r,t}f_\eps(t,\tilde  z_j)
-\cT_{r,t}f_\eps(t,0))\frac{\dif z_j}{|z_j|^{d+2/(1+2(n-j))}}\dif r\right|\dif t\\
&\leq\int^\infty_2\left|\int^s_0\!\!\!\int_{|z_j|\leq t^\gamma}|z_j|^2\|\nabla_{x_j}\nabla_{x_{j-1}}\cT_{r,t}f_\eps\|_\infty\frac{\dif z_j}{|z_j|^{d+2/(1+2(n-j))}}\dif r\right|\dif t\\
&\lesssim\int^\infty_2\left(\int^s_0(t-r)^{-2(n-j)-2}\int_{|z|\leq t^\gamma}\frac{|z_j|^2\dif z_j}{|z_j|^{d+2/(1+2(n-j))}}\dif r\right)\dif t\\
&\lesssim\int^\infty_2(t-1)^{-2(n-j)-2}t^{\gamma(2-2/(1+2(n-j)))}\dif t\lesssim 1.
\end{align*}
Combining the above calculations, we obtain the desired estimate.
\end{proof}

\bl\label{Le37}
Under  \eqref{AA}, there is a constant $C>0$ depending only on $n,\kappa, p,d$ such that for all $f\in L^\infty(\mR^{1+nd})$ with $\|f\|_\infty\leq 1$,
\begin{align}\label{EP1}
\sup_{\eps\in(0,1)}\int_{Q_1(0)}|\sP^\eps_{j2}f(s,x)-\sP^\eps_{j2}f(0,0)|^2\leq C.
\end{align}
\el
\begin{proof}
By definition, we have
\begin{align*}
&|\sP^\eps_{j2}f(s,x)-\sP^\eps_{22}f(0,0)|\leq\int^\infty_2|\e^{\lambda(s-t)}-\e^{-\lambda t}|\, \|\Delta_{x_j}^{ {1}/{(1+2(n-j))}}\cT_{s,t}f_\eps(t)\|_\infty\dif t\\
&\qquad\qquad+\int^\infty_2\e^{-\lambda t}|\Delta_{x_j}^{ {1}/{(1+2(n-1))}}\cT_{s,t}f_\eps(t,x)-\Delta_{x_j}^{ {1}/{(1+2(n-j))}}\cT_{s,t}f_\eps(t,0)|\dif t\\
&\qquad\qquad+\int^\infty_2\e^{-\lambda t}|\Delta_{x_j}^{ {1}/{(1+2(n-j))}}\cT_{s,t}f_\eps(t,0)-\Delta_{x_j}^{ {1}/{(1+2(n-j))}}\cT_{0,t}f_\eps(t,0)|\dif t\\
&\qquad\qquad=:I_1(s)+I_2(s,x)+I_3(s).
\end{align*}
Noticing that by \eqref{EK2},
\begin{align*}
\|\Delta_{x_j}^{ {1}/{(1+2(n-j))}} \cT_{s,t}f_\eps(t)\|_\infty&\lesssim (t-s)^{-1},\\
\|\nabla_{x_k}\Delta_{x_j}^{ {1}/{(1+2(n-j))}}\cT_{s,t}f_\eps(t)\|_\infty&\lesssim (t-s)^{-(2(n-k)+3)/2},
\end{align*}
we have for all $s\in[-1,1]$,
\begin{align*}
I_1&\lesssim \int^\infty_2|\e^{\lambda(s-t)}-\e^{-\lambda t}|(t-s)^{-1}\dif t\\
&\lesssim |\e^{\lambda s}-1|\int^\infty_2\e^{-\lambda t}\dif t=|\e^{\lambda s}-1|\e^{-2\lambda}/\lambda\lesssim 1,
\end{align*}
and for all $(s,x)\in Q_1(0)$,
\begin{align*}
I_2(s,x)\lesssim \int^\infty_2\sum_{k=1}^n(t-s)^{-(2(n-k)+3)/2}\dif t\lesssim 1.
\end{align*}
Moreover, by \eqref{NB6}, we have for all $s\in[-1,1]$,
$$
I_3(s)\lesssim 1.
$$
Combining the above calculations, we obtain \eqref{EP1}.
\end{proof}

Now we can give

\begin{proof}[Proof of \eqref{BN5} for $p\in(2,\infty)$]
By Lemmas \ref{Le34}, \ref{Le35} and \ref{Le37}, it is easy to see that
$\sP^\eps_j: L^\infty(\mR^{1+nd})\to BMO$ is a bounded linear operator with bound independent of $\eps$.
Estimate \eqref{BN5} for $p\in(2,\infty)$ follows by Theorem \ref{Th2} and the well-proved estimate for $p=2$.
\end{proof}

\subsection{Case: $p\in(1,2)$}We shall use duality argument to show \eqref{BN5} for $p\in(1,2)$.
Let $\cT^*_{s,t}$ be the adjoint operator of $\cT_{s,t}$, that is,
$$
\int g\cT_{s,t} f=\int f\cT^*_{s,t} g.
$$
By the definition of $\cT_{s,t}$, we have
$$
\cT^*_{s,t}f(x):=\cT^{*,a}_{s,t}f(x):=\mE f\left(\e^{(s-t)A}x-\int^t_s\e^{(s-r)A}\sigma^a_r\dif W_r\right).
$$
For $j=1,\cdots,n$, we introduce a new operator
\begin{align*}
\sQ^\eps_j f:&=\sQ^{a}_j f_\eps(s,x):=\int^t_{-\infty}\e^{\lambda(s-t)}\cT^{*,a}_{s,t}\Delta^{{1}/{(1+2(n-j))}}_{x_j} f_\eps(s,x)\dif s,
\end{align*}
where $f_\eps(t,x)=f(t,\cdot)*\varrho_\eps(x)$ so that $\sQ^\eps_j f$ is well defined for $f\in L^\infty(\mR^{1+nd})$.
Notice that $\sQ^a_j$ can be considered as the adjoint operator of $\sP^a_j$ in the sense that
$$
\int\sP^a_j f\, g=\int\sQ^a_j g\, f.
$$
As in the previous subsection, we want to show that
$$
\sQ^\eps_j\mbox{ is a bounded linear operator from $L^\infty(\mR^{1+nd})$ to $BMO$.}
$$
First of all, as in Lemma \ref{Le34} we have
$$
\fint_{Q_r(t_0,x_0)}\big|\sQ^{a}_j f_\eps(s,x)-c\big|^2=
\fint_{Q_1(0)}\big|\sQ^{\tilde a}_j\tilde f_\eps(s,x)-c\big|^2,
$$
where $\tilde a$ and $\tilde f$ are defined as in Lemma \ref{Le34}.
We aim to prove that there is a constant $C=C(n,\kappa, p,d)>0$ independent of $\eps\in(0,1)$ such that for all $f\in L^\infty(\mR^{1+nd})$
with $\|f\|_\infty\leq 1$,
$$
\fint_{Q_1(0)}\big|\sQ^{\tilde a}_j\tilde f_\eps(s,x)-c\big|^2\leq C.
$$
Below we drop $\tilde a$ and the tilde, and make the following decomposition as above,
$$
\sQ^\eps_j f(t,x)=\left(\int^t_{-2}+\int^{-2}_{-\infty}\right)\e^{\lambda(s-t)}\cT^{*}_{s,t}\Delta^{{1}/{(1+2(n-j))}}_{x_j}f_\eps(s,x)\dif s=:\sQ^\eps_{j1}f(t,x)+\sQ^\eps_{j2}f(t,x).
$$

The following lemma is crucial for treating $\sQ^\eps_{j1}$.
\bl
Let $\varphi\in C^\infty_c(\mR^{nd})$. For any $p\in[1,2]$, there are constants $C_\varphi,\beta>0$ such that
for all $h\in L^2(\mR^{nd})$ and $0<t-s\leq 3$,
\begin{align}
&\qquad\|\Delta^{{1}/{(1+2(n-j))}}_{x_j}(\cT_{s,t}(\varphi^2 h)-\varphi_{s,t}\cT_{s,t} (\varphi h))\|_p\leq C_\varphi(t-s)^{\beta-1}\|h\|_2,\label{LK0}\\
&\|\Delta^{{1}/{(1+2(n-j))}}_{x_j}(\varphi_{s,t}\cT_{s,t}(\varphi h))-\varphi_{s,t}\Delta^{{1}/{(1+2(n-j))}}_{x_j}\cT_{s,t} (\varphi h)\|_p\leq 
C_\varphi(t-s)^{\beta-1}\|h\|_2,\label{LK00}
\end{align}
where $\varphi_{s,t}(x):=\varphi(\e^{(t-s)A}x)$.
\el
\begin{proof}
(1) Let $p^{({\bf 0})}_{s,t}(y)$ be the distribution density of $X^{s,\bf 0}_t$ By the definition and the chain rule, we have
\begin{align*}
\nabla_{x_j}\cT_{s,t}f(x)&=\int_{\mR^{nd}}f(y) \nabla_{x_j}p^{({\bf 0})}_{s,t}(y-\e^{(t-s)A}\cdot)(x)\dif y\\
&=\int_{\mR^{nd}}f(y)\sum_{i=1}^j (\nabla_{y_i}p^{({\bf 0})}_{s,t})(y-\e^{(t-s)A}x)\frac{(t-s)^{j-i}}{(j-i)!}\dif y\\
&=\int_{\mR^{nd}}f(y+\e^{(t-s)A}x)\sum_{i=1}^j (\nabla_{y_i}p^{({\bf 0})}_{s,t})(y)\frac{(t-s)^{j-i}}{(j-i)!}\dif y.
\end{align*}
By this formula, we have
\begin{align}\label{HJ7}
\begin{split}
&\|\nabla_{x_j}\cT_{s,t}(\varphi^2 h)-\varphi_{s,t}\nabla_{x_j}\cT_{s,t}(\varphi h)\|_p\\
&\leq \|\nabla\varphi\|_\infty\|\varphi h\|_p\sum_{i=1}^j \int_{\mR^{nd}}|y||\nabla_{y_i}p^{({\bf 0})}_{s,t}|(y)\frac{(t-s)^{j-i}}{(j-i)!}\dif y.
\end{split}
\end{align}
By \eqref{HJ4}, we have
\begin{align*}
&\int_{\mR^{nd}}|y||\nabla_{y_i}p^{({\bf 0})}_{s,t}|(y)\dif y\lesssim(t-s)^{-(n^2d+2(n-i)+1)/2}\sum_{k=1}^n\int_{\mR^{nd}}|y_k|\e^{-c|\Theta_{(t-s)^{-1/2}}y|^2}\dif y\\
&\qquad\lesssim\sum_{k=1}^n(t-s)^{-(n^2d+2(k-i))/2}\int_{\mR^{nd}}\e^{-c|\Theta_{(t-s)^{-1/2}}y|^2/2}\dif y\lesssim\sum_{k=1}^n (t-s)^{i-k}.
\end{align*}
Substituting this into \eqref{HJ7}, we obtain
$$
\|\nabla_{x_j}\cT_{s,t}(\varphi^2 h)-\varphi_{s,t}\nabla_{x_j}\cT_{s,t}(\varphi h)\|_p\lesssim 
\|\varphi h\|_p(t-s)^{j-n}\lesssim\|h\|_2(t-s)^{j-n},
$$
and further,
$$
\|\nabla_{x_j}(\cT_{s,t}(\varphi^2 h)-\varphi_{s,t}\cT_{s,t}(\varphi h))\|_p\lesssim \|h\|_2(t-s)^{j-n}.
$$
Hence, for $j=1,2,\cdots,n-1$,
\begin{align*}
&\|\Delta^{{1}/{(1+2(n-j))}}_{x_j}(\cT_{s,t}(\varphi^2 h)-\varphi_{s,t}\cT_{s,t} (\varphi h))\|_p\\
&\leq\int_{\mR^d}\|\delta_{z_j}(\cT_{s,t}(\varphi^2 h)-\varphi_{s,t}\cT_{s,t} (\varphi h))\|_p\frac{\dif z_j}{|z_j|^{d+2/(1+2(n-j))}}\\
&\leq 2\int_{|z_j|>(t-s)^{n-j}}(\|\cT_{s,t}(\varphi^2 h)\|_p+\|\varphi_{s,t}\cT_{s,t} (\varphi h)\|_p)\frac{\dif z_j}{|z_j|^{d+2/(1+2(n-j))}}\\
&\quad+\int_{|z_j|\leq(t-s)^{n-j}}(\|\nabla_{x_j}(\cT_{s,t}(\varphi^2 h)-\varphi_{s,t}\cT_{s,t} (\varphi h))\|_p)\frac{|z_j|\dif z_j}{|z_j|^{d+2/(1+2(n-j))}}\\
&\lesssim (\|\varphi^2 h\|_p+\|\varphi h\|_p)\int_{|z_j|>(t-s)^{n-j}}\frac{\dif z_j}{|z_j|^{d+2/(1+2(n-j))}}\\
&\quad+\|h\|_2(t-s)^{j-n}\int_{|z_j|\leq(t-s)^{n-j}}\frac{|z_j|\dif z_j}{|z_j|^{d+2/(1+2(n-j))}}\lesssim \|h\|_2(t-s)^{-\frac{2(n-j)}{1+2(n-j)}}.
\end{align*}
Thus we get \eqref{LK0}. For $j=n$, \eqref{LK0} is direct by the chain rule.
\medskip\\
(2) For $j=1,\cdots,n-1$, by \eqref{EK22} we have
\begin{align*}
&\|\Delta^{{1}/{(1+2(n-j))}}_{x_j}(\varphi_{s,t}\cT_{s,t}(\varphi h))-\varphi_{s,t}\Delta^{{1}/{(1+2(n-j))}}_{x_j}\cT_{s,t} (\varphi h)\|_p\\
&\leq2\int_{\mR^d}\frac{\|\delta^{(1)}_{z_j}\varphi_{s,t}\|_\infty\|\delta^{(1)}_{z_j}\cT_{s,t}(\varphi h)\|_p}{|z_j|^{d+2/(1+2(n-j))}}\dif z_j\\
&\lesssim\|\varphi h\|_p\int_{\mR^d}\frac{1\wedge((t-s)^{n-j}|z_j|)}{|z_j|^{d+2/(1+2(n-j))}}\dif z_j\lesssim \|h\|_2(t-s)^{\frac{2(n-j)}{1+2(n-j)}},
\end{align*} 
which gives \eqref{LK00}.  For $j=n$, \eqref{LK00} is direct by the chain rule.
\end{proof}
\bl\label{Le354}
Under  \eqref{AA}, there is a constant $C=C(n,\kappa, p,d)>0$
such that for all $f\in L^\infty(\mR^{1+nd})$ with $\|f\|_\infty\leq 1$,
\begin{align}
\sup_{\eps\in(0,1)}\int_{Q_1(0)}|\sQ^{\eps}_{j1} f(s,x)|^2\leq C.\label{BN4'}
\end{align}
\el
\begin{proof}
For $t\in[-2,1]$, define
$$
u_\eps(t,x):=\sQ^\eps_{j1}f(t,x)=\int^t_{-\infty}\e^{\lambda(s-t)}\cT^*_{s,t}\Delta^{{1}/{(1+2(n-j))}}_{x_j} ((1_{[-2,1]}f_\eps)(s))(x)\dif s.
$$
Let $\varphi$ be a nonnegative smooth cutoff function in $\mR^{nd}$ with $\varphi(x)=1$ for $|x|\leq 1$ and  $\varphi(x)=0$ for $|x|>2$. We have
\begin{align*}
\|u\|_{L^2(Q_1(0))}&\leq\|1_{[-2,1]}u\varphi^2\|_{2}=\sup_{\|h\|_2\leq 1}\int_{\mR^{1+nd}}1_{[-2,1]}u\varphi^2h\\
&=\sup_{\|h\|_2\leq 1}\int_{\mR^{1+nd}}1_{[-2,1]}f_\eps\Delta^{{1}/{(1+2(n-j))}}_{x_j}
\int^1_\cdot \e^{\lambda(\cdot-t)}\cT_{\cdot,t}(\varphi^2 h (t))\dif t\\
&\leq\sup_{\|h\|_2\leq 1}\left\|1_{[-2,1]}\Delta^{{1}/{(1+2(n-j))}}_{x_j}\int^1_{\cdot}\e^{\lambda(s-t)}\cT_{s,t}(\varphi^2 h (t))\dif t\right\|_1.
\end{align*}
Since the support of $\varphi_{s,t}(x)=\varphi(\e^{(t-s)A}x)$ is contained in $\big\{x: |x|\leq n\big\}$ for $|t-s|\leq 3$, 
by \eqref{LK0} and \eqref{LK00}, we have for any $h\in L^2(\mR^{1+nd})$ with $\|h\|_2\leq 1$. 
\begin{align*}
&\left\|1_{[-2,1]}\Delta^{{1}/{(1+2(n-j))}}_{x_j}\int^1_{\cdot}\e^{\lambda(s-t)}\cT_{s,t}(\varphi^2 h(t))\dif t\right\|_1\\
&\quad\stackrel{\eqref{LK0}}{\lesssim}\left\|1_{[-2,1]}\Delta^{{1}/{(1+2(n-j))}}_{x_j}\int^1_{\cdot}\e^{\lambda(s-t)}\varphi_{s,t}\cT_{s,t}(\varphi h(t))\dif t\right\|_1+1\\
&\quad\stackrel{\eqref{LK00}}{\lesssim}\left\|1_{[-2,1]}\int^1_{\cdot}\e^{\lambda(s-t)}\varphi_{s,t}\Delta^{{1}/{(1+2(n-j))}}_{x_j}\cT_{s,t}(\varphi h(t))\dif t\right\|_1+1\\
&\quad\lesssim \left\|1_{[-2,1]}\int^1_{\cdot}\e^{\lambda(s-t)}\varphi_{s,t}\Delta^{{1}/{(1+2(n-j))}}_{x_j}\cT_{s,t}(\varphi h(t))\dif t\right\|_2+1\\
&\quad\stackrel{\eqref{LK00}}{\lesssim}\left\|1_{[-2,1]}\Delta^{{1}/{(1+2(n-j))}}_{x_j}\int^1_{\cdot}\e^{\lambda(s-t)}\varphi_{s,t}\cT_{s,t}(\varphi h(t))\dif t\right\|_2+1\\
&\quad\stackrel{\eqref{LK0}}{\lesssim}\left\|\Delta^{{1}/{(1+2(n-j))}}_{x_j}\int^\infty_{\cdot}\e^{\lambda(s-t)}\cT_{s,t}(1_{[-2,1]}\varphi^2 h (t))\dif t\right\|_2+1\\
&\quad\lesssim \|\varphi^2 h\|_2+1\lesssim 1,
\end{align*}
where the last step is due to \eqref{BN5} for $p=2$. The proof is complete.
\end{proof}

The following lemma is the same as in Lemma \ref{Le37}.
\bl\label{Le377}
Under  \eqref{AA}, there is a constant $C=C(n,\kappa, p,d)>0$ such that for all $f\in L^\infty(\mR^{1+nd})$ with $\|f\|_\infty\leq 1$,
$$
\sup_{\eps\in(0,1)}\int_{Q_1(0)}|\sQ^\eps_{j2}f(t,x)-\sQ^\eps_{j2}f(0,0)|^2\leq C.
$$
\el
\begin{proof}
By \eqref{EK22}, we have for all $t\in[-1,1]$,
\begin{align*}
&\int^{-2}_{-\infty}\big|\cT^*_{s,t}\Delta^{{1}/{(1+2(n-j))}}_{x_j}f_\eps(s,0)
-\cT^*_{s,0}\Delta^{{1}/{(1+2(n-j))}}_{x_j}f_\eps(s,0)\big|\dif s\\
&\quad\leq\int^{-2}_{-\infty}\!\int^t_0\big|\p_r\cT^*_{s,r}\Delta^{{1}/{(1+2(n-j))}}_{x_j}f_\eps(s,0)\big|\dif r\dif s\\
&\quad=\int^{-2}_{-\infty}\!\int^t_0\big|\sL^*_r\cT^*_{s,r}\Delta^{{1}/{(1+2(n-j))}}_{x_j}f_\eps(s,0)\big|\dif r\dif s\\
&\quad=\int^{-2}_{-\infty}\!\int^t_0\big|\tr(a_r\cdot\nabla^2_{x_n})\cT^*_{s,r}\Delta^{{1}/{(1+2(n-j))}}_{x_j}f_\eps(s,0)\big|\dif r\dif s\\
&\quad\leq C\int^{-2}_{-\infty}\!\int^t_0(r-s)^{-2}\dif r\dif s\leq C.
\end{align*}
Using this estimate and \eqref{EK22}, as in the proof of Lemma \ref{Le37}, we obtain the desired estimate.
\end{proof}

\begin{proof}[Proof of \eqref{BN5} for $p\in(1,2)$]
By Lemmas \ref{Le354} and \ref{Le377}, we know that
$$
\sQ^\eps_j: L^\infty(\mR^{1+2d})\to BMO\mbox{ is bounded with norm independent of $\eps$.}
$$
Moreover, by duality, we also have
$$
\sQ^\eps_j: L^2(\mR^{1+2d})\to L^2(\mR^{1+2d})\mbox{ is bounded with norm independent of $\eps$.}
$$
Hence, for $q=p/(p-1)\in(2,\infty)$, by Theorem \ref{Th2}, we have for some $C>0$ independent of $\eps$,
$$
\|\sQ^\eps_j f\|_q=\left\|\int^t_{-\infty}\e^{\lambda(s-t)}\cT^*_{s,t}\Delta^{{1}/{(1+2(n-j))}}_{x_j} f_\eps\dif s\right\|_q\leq C\|f\|_q.
$$
Now for $p\in(1,2)$, by Fatou's lemma, we get
\begin{align*}
\|\sP_jf\|_p&\leq \|f\|_p\sup_{\|h\|_q\leq 1}\left\|\int^t_{-\infty}\e^{\lambda(s-t)}\cT^*_{s,t}\Delta^{{1}/{(1+2(n-j))}}_{x_j} h\dif s\right\|_q\\
&\leq \|f\|_p\sup_{\|h\|_q\leq 1}\varliminf_{\eps\to 0}
\left\|\int^t_{-\infty}\e^{\lambda(s-t)}\cT^*_{s,t}\Delta^{{1}/{(1+2(n-j))}}_{x_j} h_\eps\dif s\right\|_q\leq C\|f\|_p,
\end{align*}
which gives \eqref{BN5} for $p\in(1,2)$.
\end{proof}

\end{document}